\documentclass[MSNbibl,number,citesort,seceqn,dvips]{arxbj}
\usepackage{graphicx}
\aid{0}
\volume{20}
\issue{1}
\pubyear{2014}
\firstpage{164}
\lastpage{189}
\doi{10.3150/12-BEJ480} 

\makeatletter
\newtheorem{teo}{Theorem}[section]
\newtheorem{Lem}{Lemma}[section]
\newproclaim{assum}{Assumption}
\newtheorem{Corol}{Corollary}[section]
\newremark{rem}{Remark}
\renewcommand{\widetilde}{\tilde}
\newcommand{\sgn}{\operatorname{sgn}}
\newcommand{\tr}{\operatorname{tr}}
\newcommand{\Argmin}{\operatorname{Argmin}}
\newcommand{\diag}{\operatorname{diag}}
\newcommand{\DIFadd}[1]{#1} 
\makeatother

\begin{document}
\begin{frontmatter}

\title{Estimating spatial quantile regression with functional
coefficients: A robust semiparametric framework}
\runtitle{Spatial quantile regression with functional
coefficients}

\begin{aug}
\author[1]{\fnms{Zudi} \snm{Lu}\corref{}\thanksref{1}\ead[label=e1]{zudi.lu@adelaide.edu.au}\ead[label=u1,url]{http://www.adelaide.edu.au/directory/Zudi.Lu}},%
\author[2]{\fnms{Qingguo} \snm{Tang}\thanksref{2,e2}\ead[label=e2,mark]{tangqig@yahoo.com.cn}} \and
\author[2]{\fnms{Longsheng} \snm{Cheng}\thanksref{2,e3}\ead[label=e3,mark]{cheng\_longsheng@163.com}}
\runauthor{Z. Lu, Q. Tang and L. Cheng} 
\address[1]{School of Mathematical Sciences, The University of
Adelaide, SA5005,
Australia.\\ \printead{e1,u1}}
\address[2]{School of Economics and Management, Nanjing University
of Science and Technology, Nanjing 210094, China. \printead{e2,e3}}
\end{aug}

\received{\smonth{1} \syear{2012}}
\revised{\smonth{9} \syear{2012}}

%
\begin{abstract}
This paper considers an estimation of semiparametric functional
(varying)-coefficient quantile regression
with spatial data. A general robust framework is developed
that treats quantile regression for spatial data in a natural
semiparametric way. The local M-estimators of the unknown
functional-coefficient functions are proposed by using local
linear approximation, and their asymptotic
distributions
are then established under weak spatial mixing conditions allowing
the data processes to be either stationary or nonstationary with
spatial trends.
Application to a soil data set is demonstrated with interesting
findings that go beyond traditional analysis.
\end{abstract}

%

\end{frontmatter}

\section{Introduction}

Spatial data, which are collected at different sites on the surface of
the earth, arise in various areas of research, including econometrics,
epidemiology, environmental
science, image analysis, oceanography and many others. Numerous
applications of spatial models and important developments in the
general area of spatial statistics under linear correlation structures
can be found in \cite{s5,s1,s2,s12}, and a more recent
comprehensive review by Gelfand \textit{et al.}
\cite{s11b}, among others. However, linear correlation structures may not
be always reasonable in spatial applications.
In the last ten years, efforts have been made in the literature to
explore nonlinear relationship in spatial data. See,
for example, \cite{s28,s29,s14,s11,s31,s30}, who explored the nonlinear spatial
interdependence from the perspective of conditional mean
regressions. Differently from these references, Hallin \textit{et al.} \cite{s15} recently proposed to investigate the nonlinear
spatial interaction by using conditional quantile regression, showing
that spatial quantile regression can provide much
more information on spatial data than the conditional mean regression analysis.
In this paper, following the above efforts, we aim to develop a
structure of spatial quantile regression allowing
functional coefficients, \textit{under a robust semiparametric
framework}, to reduce the ``curse of dimensionality''
that spatial quantile regression analysis suffers from when the
dimension of the covariates is higher than 3. We will
demonstrate in Section~\ref{sec5} that the proposed semiparametric
functional-coefficient spatial quantile structure will be
useful in the analysis of a soil data set.

To make our results widely applicable, we shall consider the quantile
regression for
spatial data in a general context. Firstly, we treat data as observed
over a space of general dimension $N$.
Denote the set of integer lattice points in $N$-dimensional Euclidean space
by $Z^{N}$, where $N\geq1$ and $Z=\{0,\pm1, \pm2, \ldots\}$.
A point $\mathbf{i}=(i_{1},\ldots,i_{N})$ in $Z^{N}$ is referred to as
a site.
Spatial data are modeled as finite realizations of vector stochastic processes
indexed by $\mathbf{i}\in Z^{N}$, that is, random fields. We
will consider strictly stationary $(d+k+1)$-dimensional random fields
of the form
\[
\bigl\{(Y_{\mathbf{i}},X_{\mathbf{i}},U_{\mathbf{i}})\dvtx  \mathbf{i}\in
Z^{N}\bigr\},
\]
where $Y_{\mathbf{i}}$, with values in $R$, $X_{\mathbf{i}}$,
with values in $R^{d}$, and $U_{\mathbf{i}}$, with values in $R^{k}$,
are defined over a probability space
$(\Omega,\mathcal{F},P)$.

Secondly, we treat spatial quantile regression in a general context of
robust spatial regression. In a number of applications,
a crucial problem consists in describing and analyzing the influence
of the covariates $(U_{\mathbf{i}},X_{\mathbf{i}})$ on the
real-valued response $Y_{\mathbf{i}}$. In spatial context, this
study is particularly difficult due to the possibly highly
complex spatial dependence among the various sites. 
The traditional approach to this problem consists in assuming that
$Y_{\mathbf{i}}$ has finite expectation, so that spatial
conditional mean regression function $g\dvtx  (x,u)\mapsto g(x,u):=
E[Y_{\mathbf{i}}|X_{\mathbf{i}} = x,U_{\mathbf{i}}=u]$ may be well
defined and clearly carries relevant information on the dependence
of $Y$ on $X$ and $U$ (cf., \cite{s28,s29,s14}). Differently, Hallin \textit{et al}. \cite{s15} proposed
spatial conditional quantile regression, defined by
%
\begin{equation}\label{qtau}
 q_{\tau}\dvtx  (x,u)\mapsto q_{\tau}(x,u):= Q[Y_{\mathbf
{i}}|X_{\mathbf{i}} =
x,U_{\mathbf{i}}=u],
\end{equation}
which provides more comprehensive information
on the dependence of $Y$ on $X$ and $U$ through different $0<\tau<1$
(see \cite{s25} and \cite{s44}),
where $q_{\tau}(x,u)$ satisfies $P[Y_{\mathbf{i}}<q_{\tau
}(x,u)|X_{\mathbf{i}} =
x,U_{\mathbf{i}}=u]=\tau$; see also the robust spatial conditional
regression in \cite{s26}. As is \DIFadd{well known} in the
nonparametric literature, when $d+k> 3$, both spatial regression
functions $g(x,u)$ and $q_{\tau}(x,u)$ can not be well estimated
nonparametrically with reasonable accuracy owing to the curse of
dimensionality. Because of complex spatial interaction, this issue
on how to avoid the curse of dimensionality becomes particularly
important, which has been addressed by Gao \textit{et al}. \cite{s11}
and Lu \textit{et al}. \cite{s31} for spatial conditional mean regression
$g(x,u)$ under \textit{least squares partially linear and additive
approximation structures}, respectively. In this paper, we are
particularly concerned with avoiding the curse of dimensionality
for spatial quantile regression analysis, and, for generality,
consider a general spatial regression that takes conditional
quantile regression
$Q(Y_{\mathbf{i}}|U_{\mathbf{i}},X_{\mathbf{i}})$ as a special case,
to be
approximated by a popular linear
structure allowing for functional coefficients in the form
%
\begin{equation}\label{1.2a}
\Psi(U_{\mathbf{i}},X_{\mathbf{i}})
=X_{\mathbf{i}1}\beta_{1}(U_{\mathbf{i}})+\cdots+X_{\mathbf{i}d}\beta
_{d}(U_{\mathbf{i}}),
\end{equation}
with the functional coefficients $\beta_{j}(\cdot)$'s defined by minimizing
%
\begin{equation}\label{1.2b}
 E\rho\bigl(Y_{\mathbf{i}}-\Psi(U_{\mathbf{i}},X_{\mathbf
{i}})\bigr)=E\rho\bigl(Y_{\mathbf{i}}
-X_{\mathbf{i}1}\beta_{1}(U_{\mathbf{i}})-\cdots-X_{\mathbf{i}d}\beta
_{d}(U_{\mathbf{i}})\bigr),
\end{equation}
associated with $\rho(y)$ by which we denote hereafter for a
general loss function [see Section~\ref{sec2}], over a
class of functional coefficient linear functions of the form
$\Psi(U_{\mathbf{i}},X_{\mathbf{i}})$ in (\ref{1.2a}). In
the subsequent, when considering $\tau$th quantile regression, we
will denote by $\rho_{\tau}(z)=|z|+(2\tau-1)z$
with
$0<\tau<1$, instead of $\rho(\cdot)$, for the loss function, under which
the resulting $\Psi(U_{\mathbf{i}},X_{\mathbf{i}})$ in
(\ref{1.2a}) is the spatial quantile regression with
functional coefficients that we are mainly concerned with in this
paper.
Let $X_{\mathbf{i}}=(X_{\mathbf{i}1},\ldots,X_{\mathbf{i}d})^{T}$. As
in traditional
linear regression when a baseline effect is desired, we set $X_{i1}
\equiv1$.
The regime $U_{\mathbf{i}}$ is a vector of explanatory variables,
and $\beta_{1}(u)$, $\ldots$\,, $\beta_{d}(u)$ are unknown smooth
functions of $u$ to be estimated, with the dimension $k$ of
$U_{\mathbf{i}}$ usually small, say $k=1$ or $2$.

Functional (varying)-coefficient regressions are a useful
extension of the classical linear regressions. One of the
advantages of such models is that the effects of the regressor
vector $X_\mathbf{i}$ can be well measured by the functional
coefficients through $U_\mathbf{i}$ and the dimensionality curse
is therefore reduced when $k$ is small. Functional coefficient
regression models are popular in traditional regression and time
series analysis.
A comprehensive theory in the nonspatial case has been well explored,
see, for example, \cite{s21,s42,s9,s4,s22,s43,38,s10,36}. However,
the varying coefficient models with spatial data are still rather
rarely investigated in the literature. Some exceptions include
an extension of the useful semiparametric model
studied by Moyeed and Diggle \cite{s32}, where the intercept coefficient
$\beta_{1}$ is
assumed to be time-varying, while $\beta_{2},\ldots,\beta_{d}$ are
constants; see also \cite{s30} for a
varying-coefficient spatiotemporal model under the least squares mean
regression perspective.

In this paper, we will develop in Section~\ref{sec2} \textit{a general robust
$M$-type semiparametric framework for approximating a
spatial conditional regression}, under $\rho(\cdot)$, \textit{by the
linear
structure with functional coefficients, $\Psi(u,x)$ in (\ref{1.2a})},
via minimizing (\ref{1.2b}). We apply local linear method to
approximate the unknown coefficient functions $\beta_{r}(u),r=1,\ldots
,d$ and
obtain their local M-estimators in Section~\ref{sec2.1}. The main results on asymptotic
distribution for the local M-estimators of $\beta_r(u)$'s at both
interior and boundary
points with stationary spatial data are established in Section~\ref{sec2.2}.
Applications of the main results to conditional quantile
coefficient functions and robust conditional regression
coefficient functions will be presented in Section~\ref{sec3}.
Section~\ref{sec4}
extends
the main
results to the case of allowing a nonstationary random field with
spatial trend, which is of importance in practice.
A real data example will be reported in Section~\ref{sec5}.
The proofs of the main theorems are relegated in \hyperref[Aproofs]{Appendix}, with
details of the proof of necessary lemmas provided in the
supplementary material \cite{s31a}.

\section{Spatial quantile regression under general M-estimation
framework: Asymptotic results}\label{sec2}
Consider a rectangular sampling region by
\[
G_{N}=\bigl\{\mathbf{i}=(i_{1},\ldots,i_{N})
\in Z^{N}\dvtx  1\leq i_{l}\leq n_{l}, l=1,\ldots,N\bigr\},
\]
with $\mathbf{n}=(n_{1},\ldots,n_{N})$. In this paper, we write
$\mathbf{n}\rightarrow\infty$ if $n_{l}\rightarrow\infty$ for
some $1\leq l\leq N$. 
Assume that we observe $(Y_{\mathbf{i}}, X_{\mathbf{i}},U_{\mathbf
{i}})$ on $G_{N}$.
The total sample size is thus $\tilde{\mathbf{n}}=\prod_{l=1}^{N}n_{l}$.
We will assume that
$\{(Y_{\mathbf{i}},X_{\mathbf{i}},U_{\mathbf{i}})\}$ satisfies the following
mixing condition as defined in the literature (cf., \cite{s39,s14,s15}).
Let $S$ and $S'$ be two sets of sites.
The Borel fields $\mathcal{B(S)}=\mathcal{B}((Y_{\mathbf{i}},X_{\mathbf{i}},
U_{\mathbf{i}})\dvtx  \mathbf{i}\in S)$ and $\mathcal{B(S')}=\mathcal
{B}((Y_{\mathbf{i}},
X_{\mathbf{i}},U_{\mathbf{i}})\dvtx  \mathbf{i}\in S')$ are the $\sigma
$-fields generated by
$(Y_{\mathbf{i}},X_{\mathbf{i}},U_{\mathbf{i}})$ with $\mathbf{i}$
being the elements of $S$ and $S'$, respectively. Let $d(S,S')$ be the Euclidean
distance between $S$ and $S'$. Then the spatial mixing defines that
there exists a function $\varphi(t)\downarrow0$ as $t\rightarrow
\infty$, such that whenever $S, S'\subset Z^{N}$,
%
\begin{eqnarray}\label{1.1}
\alpha\bigl(\mathcal{B}(S),\mathcal{B}\bigl(S'\bigr)\bigr)
& =&\sup\bigl\{\bigl|P(AB)-P(A)P(B)\bigr|,
A\in\mathcal{B}(S),B\in\mathcal{B}\bigl(S'\bigr)\bigr\} \nonumber\\
[-8pt]\\[-8pt]
&\leq&\chi\bigl(\operatorname{Card}(S),\operatorname{Card}\bigl(S'\bigr)\bigr)\varphi\bigl(d\bigl(S,S'\bigr)\bigr),\nonumber
\end{eqnarray}
where $\operatorname{Card}(S)$ denotes the cardinality of $S$,
and $\chi$ is a symmetric positive function nondecreasing in each
variable. If $\chi\equiv1$, then
$\{(Y_{\mathbf{i}},X_{\mathbf{i}},U_{\mathbf{i}})\}$ is called
strongly mixing.

\subsection{A general M-type semiparametric framework}\label{sec2.1}
Consider $u\in E^{0}=\{u=(u_{1},\ldots, u_{k})\mid u_{\ast j}\leq
u_{j}\leq u_{j}^{\ast},1\leq j\leq k\}$,
where $u_{\ast j}$ and $u^{\ast}_{ j}$ are\vspace*{1pt} constants of lower and
upper limits of \DIFadd{$u_j,$} respectively.
Let 
$\beta_{r}(u)$, $r=1,\ldots,d$, in (\ref{1.2a}) be defined by minimizing
(\ref{1.2b}) with $\rho(\cdot)$. Then, given $u_{0}\in E^{0}$, for $u$
in the neighborhood of $u_{0}$, we can use $
a_{r}+b_{r}^{T}(u-u_{0})$
to approximate the unknown coefficient function $\beta
_{r}(u)$ ($r=1,\ldots,d$), where $b_{r}=(b_{r1},\ldots,b_{rk})^{T}$.
Based on spatial observations $\{(Y_{\mathbf{i}},X_{\mathbf
{i}},U_{\mathbf{i}})\dvtx  \mathbf{i}\in
G_{N}\} $, by using
the idea of local linear fitting (see, e.g., \cite{s6} and \cite{35}),
we solve the following minimization problem 
%
\begin{equation}\label{2.1}
\min_{a_{r},b_{r},r=1,\ldots,d}\sum_{\mathbf{i}\in
G_{N}}\rho\left(Y_{\mathbf{i}}-\sum_{r=1}^{d}\bigl[a_{r}+(U_{\mathbf
{i}}-u_{0})^{T}b_{r}\bigr]X_{\mathbf{i}r}\right)
K\left(\frac{U_{\mathbf{i}}-u_{0}}{h_{\mathbf{n}}}\right),
\end{equation}
where $K(\cdot)$ is a given kernel function, $h_{\mathbf{n}}$ is a
chosen bandwidth. Let
$\hat{a}_{r},\hat{b}_{r}$, $r=1,2,\ldots,d$, be the minimizer of (\ref
{2.1}). 
Then the M-estimator $\hat{\beta}(u_{0})$ of
$\beta(u_{0})=(\beta_{1}(u_{0}),\ldots,\beta_{d}(u_{0}))^{T}$, which
minimizes (\ref{1.2b}) for 
$\rho(\cdot)$, is defined by
%
\begin{equation}\label{2.2}
\hat{\beta}(u_{0})
=\bigl(\hat{\beta}_{1}(u_{0}),\ldots,\hat{\beta}_{d}(u_{0})\bigr)^{T}=\hat
{a}=(\hat{a}_{1},\ldots,\hat{a}_{d})^{T}.
\end{equation}

Typical choices for $\rho$ are convex and symmetric about 0. Here\DIFadd
{, } we only require $\rho$ a convex function so that the optimisations
(\ref{2.1})
are well defined and the problem of local
minima is avoided. It can be asymmetric. For example, 
an estimator with
$\rho_{\tau} (z)$ for $0<\tau<1$ gives the $\tau$th conditional
quantile of $Y$, defined in (\ref{qtau}).
For robustness consideration, we may take $\rho$ having a bounded
derivative $\rho'(z)=\max\{-1,\min\{z/c,1\}\}$, 
$c>0$; see \cite{s23} and \cite{s16} for more details
about the robustness of M-estimators.

\subsection{Asymptotic results}\label{sec2.2}
 In this subsection, we state asymptotic properties of the
estimates $\hat{\beta}(u_{0})$. Let $\psi(z)
$ be the derivative function of $\rho(z)$ with respect to $z$ almost everywhere.
The following assumptions are
required for our asymptotic results.
%
\begin{assum}\label{asum1}
 The random field $\{(Y_{\mathbf{i}},X_{\mathbf
{i}},U_{\mathbf{i}})\dvtx  \mathbf{i}\in
Z^{N}\}$ is strictly stationary. For all distinct $\mathbf{i}$ and
$\mathbf{j}$
in $Z^{N}$, the random variables $U_{\mathbf{i}}$ and $U_{\mathbf{j}}$
admit a joint
density 
$f_{\mathbf{i},\mathbf{j}}(u,v)\leq C_{0}$ uniformly with respect to
$\mathbf{i},\mathbf{j}\in Z^{N}$ and 
$u,v\in E^{0}$,
where $C_{0}$ is some positive constant. The marginal density $f(u)$
of $U_{\mathbf{i}}$ is continuous and bounded away from 0 uniformly
over $E^{0}$.
\end{assum}
%
\begin{assum}\label{asum2}
All functions, $\beta_{r}(u)$'s, are twice
continuously differentiable
in a neighborhood of $u_{0}$, for $r=1,\ldots,d$.
\end{assum}
%
\begin{assum}\label{asum3}
 The convex loss function $\rho(\cdot)$
satisfies, for some $\delta>0$,
\[
E\psi(\varepsilon_{\mathbf{i}}|X_{\mathbf{i}}, U_{\mathbf{i}})=0,\qquad
E\bigl(\bigl|\psi(\varepsilon_{\mathbf{i}})\bigr|^{2+\delta}|X_{\mathbf{i}},
U_{\mathbf{i}}\bigr)\leq C_{1},
\]
where
$\varepsilon_{\mathbf{i}}=Y_{\mathbf{i}}-X_{\mathbf{i}}^{T}\beta
(U_{\mathbf{i}})$
and $C_{1}>0$ is a constant. Furthermore, there exist some
function $\phi(\cdot)$ and constant $\bar{c}_{1}>0$ such that
$|E(\psi(\varepsilon_{\mathbf{i}}+z)|X_{\mathbf{i}},
U_{\mathbf{i}})-\phi(X_{\mathbf{i}}, U_{\mathbf{i}})z|\leq C_{1}z^{2}$
for any $|z|\leq\bar{c}_{1}$.
\end{assum}
%
\begin{assum}\label{asum4}
 There exist constants $0<\bar{c}_{2}, C_{2}<\infty
$ such that
\[
\bigl|\psi(v+z)-\psi(v)\bigr|\leq C_{2}, \qquad E\bigl(\bigl[\psi(\varepsilon_{\mathbf{i}}+z)-\psi
(\varepsilon_{\mathbf{i}})\bigr]^{2}|X_{\mathbf{i}},
U_{\mathbf{i}}\bigr)\leq C_{2}|z|
\]
for any $|z|\leq\bar{c}_{2}$ and
$v\in R^{1}$.
\end{assum}
%
\begin{assum}\label{asum5}
 The bandwidth $h_{\mathbf{n}}$ satisfies that
$h_{\mathbf{n}}\leq C_{3}\tilde{\mathbf{n}}^{-1/(k+4)}$ for some
positive constant $C_{3}$ and
$\tilde{\mathbf{n}}h_{\mathbf{n}}^{k}\rightarrow+\infty$ as
$\mathbf{n}\rightarrow\infty$.
\end{assum}
%
\begin{assum}\label{asum6}
$\max_{\mathbf{i}\in
G_{N}}\|X_{\mathbf{i}}\|=\mathrm{o}_{p}((\tilde{\mathbf{n}}h_{\mathbf{n}}^{k})^{1/2})$,
$\max_{\mathbf{i}\in
G_{N}}\|X_{\mathbf{i}}\|=\mathrm{o}_{p}(h_{\mathbf{n}}^{-2})$\break  and
$E\|X_{\mathbf{i}}\|^{4+2\delta}
<\infty$. 
\end{assum}
%
\begin{assum}\label{asum7}
 The kernel function $K(\cdot)\geq0$
is a bounded symmetric function with a compact support $\tilde
{M}=[-M_{1},M_{1}]\times\cdots\times[-M_{k},M_{k}]$
and $\int_{\tilde{M}}uu^{T}K(u)\,\mathrm{d}u$ is positive definite.
\end{assum}
%
\begin{assum}\label{asum8}
The function $\chi(\cdot, \cdot)$ and $\varphi$
satisfy that
$\chi(n',n'')\leq\min(n',n'')$ and
%
\begin{equation}\label{mix}
\lim_{k\rightarrow\infty}k^{a}\sum_{z=k}^{\infty}z^{N-1}\bigl\{\varphi(z)\bigr\}
^{\delta/(2+\delta)}
=0
\end{equation}
for some constant
$a>(4+\delta)N/(2+\delta)$.
\end{assum}
%
\begin{assum}\label{asum9}
 $\min_{1\leq l\leq N}n_{l}\rightarrow\infty$
and there exist two sequences of positive integer vectors,
$\mathbf{p}=\mathbf{p}_{\mathbf{n}}=(p_{1},\ldots, p_{N})\in
Z^{N}$ and $\mathbf{q}=\mathbf{q}_{\mathbf{n}}=(q,\ldots,q)\in
Z^{N}$, with $q\rightarrow\infty$ such that $q/p_{l}\rightarrow
0$ and $n_{l}/p_{l}\rightarrow\infty$ for all $l=1,\ldots,N$, and
$\tilde{\mathbf{p}}=\prod_{l=1}^{N}p_{l}=\mathrm{o}((\tilde{\mathbf{n}}h_{\mathbf
{n}}^{k})^{1/2})$,
$\tilde{\mathbf{n}}\varphi(q)\rightarrow0$. Furthermore,
$qh_{\mathbf{n}}^{\delta k/[a(2+\delta)]}>1$.
\end{assum}

The above assumptions are standard in the setting of local
smoothers needed for asymptotics. See \cite{s14}, for
example, for Assumptions \ref{asum1}, \ref{asum2}, \ref{asum7}, \ref{asum8} and \ref{asum9} in the spatial
context. Assumptions \ref{asum3} and~\ref{asum4} are easily checked if
the score function $\psi$ is differentiable, but they cover
nondifferentiable case including the least absolute
deviation estimator with $\psi(z)=\sgn(z)$.
Assumptions \ref{asum5} and~\ref{asum6} can be found
in \cite{s9}, where the moment condition on $X_{\mathbf
{i}}$ is technical
for the establishment of asymptotic properties
in varying coefficient setting. The bounded support
restriction on $K(\cdot)$ is technical
and can be relaxed by using such kernels with light tails as Gaussian
kernel.

To state our main results, we let
\begin{eqnarray*}
\Phi(u)&=&E\bigl(\phi(X_{\mathbf{i}},U_{\mathbf{i}})X_{\mathbf{i}}X_{\mathbf
{i}}^{T}|U_{\mathbf{i}}=u\bigr),\qquad
\Sigma(u)=E\bigl(E\bigl(\psi^{2}(\varepsilon_{\mathbf{i}})|X_{\mathbf
{i}},U_{\mathbf{i}}\bigr)
X_{\mathbf{i}}X_{\mathbf{i}}^{T}|U_{\mathbf{i}}=u\bigr),
\\
\zeta(u_{0})&=&\bigl(\zeta_{1}(u_{0}),\ldots,\zeta_{d}(u_{0})\bigr)^{T}, \qquad \zeta
_{r}(u_{0})=\tr\left(\ddot{\beta}_{r}(u_{0})\int_{\tilde
{M}}uu^{T}K(u)\,\mathrm{d}u\right),\qquad
r=1,\ldots,d,
\end{eqnarray*}
where $\ddot{\beta}_{r}(u_{0})$ is the second derivative of
$\beta_{r}(u)$ at $u=u_0$.
%
\begin{teo}\label{thm2.1}
Assume that Assumptions \ref{asum1}--\ref{asum9} hold and
$\Phi(u), \Sigma(u)$ are continuous in some neighborhood of
$u_{0}$ and $\Phi(u_{0})$ is positive definite. If $u_{0}$ is an
interior point of the support of the design density $f(u)$, then,
as $\mathbf{n}\rightarrow\infty$,
\[
\sqrt{\tilde{\mathbf{n}}h_{\mathbf{n}}^{k}}\left(\hat{\beta
}(u_{0})-\beta(u_{0})-
\frac{h_{\mathbf{n}}^{2}}{2\mu_{0}}\zeta(u_{0})\right)\rightarrow_{d}
N\biggl(0,
\frac{\nu_{0}}{f(u_{0})\mu_{0}^{2}}\Phi^{-1}(u_{0})\Sigma(u_{0})\Phi
^{-1}(u_{0})\biggr),
\]
where $\mu_{0}=\int_{\tilde{M}}K(u)\,\mathrm{d}u$,
$\nu_{0}=\int_{\tilde{M}}K^{2}(u)\,\mathrm{d}u$, and $\rightarrow_{d}$ means
convergence in distribution.
\end{teo}

Theorem~\ref{thm2.1} gives the asymptotic distribution of
the estimator of $\beta(u_{0})$ at an interior point.
Next\DIFadd{,} we study the asymptotic behavior of the estimator at the
boundary of the support $E^{0}$ of $f(u)$. Suppose
$u_{*}=(u_{*1},\ldots,u_{*k})^{T}$ is a boundary point. Take
$u_{h}=u_{*}+ch_{\mathbf{n}}$, where $c=(c_{1},\ldots,c_{k})^{T}$
satisfies that $0\leq c_{l}<M_{l},l=1,\ldots,k$. 
Let $\bar{M}=[-c_{1},M_{1}]\times\cdots\times[-c_{k},M_{k}]$, $
\bar{\zeta}(u_{*})=(\bar{\zeta}_{1}(u_{*}),\ldots,\bar{\zeta}_{d}(u_{*}))^{T},$
$\bar{\zeta}_{r}(u_{*})=\tr(\ddot{\beta}_{r}(u_{*})\int_{\bar{M}}vv^{T}K(v)\,\mathrm{d}v),
r=1,\ldots,d, $
$\bar{\zeta}^{^{(l)}}(u_{*})=(\bar{\zeta}_{1}^{(l)}(u_{*}),\ldots,\bar
{\zeta}_{d}^{(l)}(u_{*}))^{T},
$
$\bar{\zeta}_{r}^{(l)}(u_{*})=\tr(\ddot{\beta}_{r}(u_{*})\int_{\bar
{M}}v_{l}vv^{T}K(v)\,\mathrm{d}v),
l=1,\ldots,k, $ and
\[
\begin{array}{l}
\Delta_{c}=\left(
\begin{array}{cc}
\displaystyle \int_{\bar{M}}K(u)\,\mathrm{d}u & \displaystyle \int_{\bar{M}}u^{T}K(u)\,\mathrm{d}u\\\noalign{\vspace*{2pt}}
\displaystyle \int_{\bar{M}}uK(u)\,\mathrm{d}u & \displaystyle \int_{\bar{M}}uu^{T}K(u)\,\mathrm{d}u
\end{array}
\right),\qquad
\bar{\Delta}_{c}=\left(
\begin{array}{cc}
\displaystyle \int_{\bar{M}}K^{2}(u)\,\mathrm{d}u & \displaystyle \int_{\bar{M}}u^{T}K^{2}(u)\,\mathrm{d}u\\\noalign{\vspace*{2pt}}
\displaystyle \int_{\bar{M}}uK^{2}(u)\,\mathrm{d}u & \displaystyle \int_{\bar{M}}uu^{T}K^{2}(u)\,\mathrm{d}u
\end{array}
\right).
\end{array}
\]

\begin{teo}\label{thm2.2}
Assume that Assumptions \ref{asum1}--\ref{asum9} hold in
some right neighborhood of $u_{*}$ and $\Phi(u), \Sigma(u)$ are
continuous in some right neighborhood of $u_{*}$ and $\Phi(u_{*})$
is positive definite. Suppose $\Delta_{c}$ is invertible. Then, as
$\mathbf{n}\rightarrow\infty$,
\begin{eqnarray*}
&&\sqrt{\tilde{\mathbf{n}}h_{\mathbf{n}}^{k}}\left(\hat{\beta
}(u_{h})-\beta(u_{h})
-\frac{h_{\mathbf{n}}^{2}}{2}\bigl[\delta_{11}
\bar{\zeta}(u_{*})+\sum_{l=1}^{k}\delta_{1(l+1)}\bar{\zeta
}^{^{(l)}}(u_{*})\bigr]\right)
\\
&&\quad \rightarrow_{d} N\left(0,
\frac{\lambda_{11}}{f(u_{*})}\Phi^{-1}(u_{*})\sum(u_{*})\Phi
^{-1}(u_{*})\right),
\end{eqnarray*}
where $\delta_{ij}$ denotes the $(i,j)$th entry of $
\Delta_{c}^{-1}$ and $\lambda_{11}$ denotes the $(1,1)$th entry of
$\Delta_{c}^{-1}\bar{\Delta}_{c}\Delta_{c}^{-1}$.
\end{teo}

The proofs of Theorems \ref{thm2.1} and \ref{thm2.2} are postponed
to \hyperref[appendix]{Appendix}.

Theorem~\ref{thm2.2} shows that for local linear estimator
the convergence rate at the points near the boundary is the same as
that for interior points.
Hence for local linear estimator near the boundary no adjustments
\DIFadd{are} required.

For the mixing coefficient $\varphi(t)$, if it decays at an
algebraic rate, that is, $\varphi(t)=\mathrm{O}(t^{-\mu})$ for some
$\mu>2(3+\delta)N/\delta$, we can choose constant $a$ such that
$(4+\delta)N/(2+\delta)<a<\mu\delta/(2+\delta)-N$, then, as
$l\rightarrow\infty$, it holds that
\begin{eqnarray*}
l^{a}\sum_{z=l}^{\infty}z^{N-1}\bigl\{\varphi(z)\bigr\}^{\delta/(2+\delta)}
&\leq& Cl^{a}\sum_{z=l}^{\infty}z^{N-1}z^{-\mu\delta/(2+\delta)} \\
 &\leq&
Cl^{a}l^{N-\mu\delta/(2+\delta)}=Cl^{-[\mu\delta/(2+\delta
)-a-N]}\rightarrow
0,
\end{eqnarray*}
and so (\ref{mix}) holds. Using the similar arguments to those
used in the proof of Theorem~3.3 of \cite{s14},
Assumption~\ref{asum9} can be much simplified, and we have the following
\DIFadd{corollary.} 
%
\begin{Corol}\label{cor1}
 Assume that Assumptions \ref{asum1}--\ref{asum7} hold and
$\Phi(u), \Sigma(u)$ are continuous in some neighborhood of
$u_{0}$ and $\Phi(u_{0})$ is positive definite. Suppose
$\Delta_{c}$ is invertible and $\chi(n',n'')\leq\min(n',n'')$ and
$\varphi(t)=\mathrm{O}(t^{-\mu})$ for some $\mu>2(3+\delta)N/\delta$. Let
the sequence of positive integers $q
=q_{\mathbf{n}}\rightarrow\infty$ and the bandwidth
$h_{\mathbf{n}}$ such that
$\tilde{\mathbf{n}}q^{-\mu}\rightarrow0$, $q=\mathrm{o}(\min_{1\leq i\leq
N}(n_{i}h_{\mathbf{n}}^{k/N})^{1/2})$ and $qh_{\mathbf{n}}^{\delta
k/[a(2+\delta)]}>1$. Then, as $\mathbf{n}\rightarrow\infty$, the
conclusions of Theorems \ref{thm2.1} and \ref{thm2.2} still hold.
\end{Corol}

Further, if the mixing coefficient decays at a geometric rate,
that is, $\varphi(t)=\mathrm{O}(\mathrm{e}^{-\nu t})$ for some $\nu>0$, then
similarly to 
Theorem~3.4 of \cite{s14}, Assumptions \ref{asum8} and \ref{asum9} can also
be simplified and we have the
following \DIFadd{corollary.} 
%
\begin{Corol}\label{cor2}
Assume that Assumptions \ref{asum1}--\ref{asum7} hold and
$\Phi(u), \Sigma(u)$ are continuous in some neighborhood of
$u_{0}$ and $\Phi(u_{0})$ is positive definite. Suppose
$\Delta_{c}$ is invertible and $\chi(n',n'')\leq\min(n',n'')$ and
$\varphi(t)=\mathrm{O}(\mathrm{e}^{-\nu t})$ for some $\nu>0$. If
\[
\min_{1\leq i\leq
N}\bigl\{\bigl(n_{i}h_{\mathbf{n}}^{k/N}\bigr)^{1/2}\bigr\}h_{\mathbf{n}}^{\delta
k/[a(2+\delta)]}(\ln\tilde{\mathbf{n}})^{-1} \rightarrow\infty
\]
for some constant $a>(4+\delta)N/(2+\delta)$, then, as
$\mathbf{n}\rightarrow\infty$, the conclusions of
Theorems \ref{thm2.1} and~\ref{thm2.2} still hold.
\end{Corol}
%
\begin{rem}\label{rem1}
 Another way for $\mathbf{n}$ to tend to infinity
is the so called isotropic one, where all components of
$\mathbf{n}$ tend to infinity at the same rate. We write
$\mathbf{n}\Rightarrow\infty$ if $\mathbf{n}\rightarrow\infty$ and
$|n_{j}/n_{l}|<C_{4}$ for some $0<C_{4}<\infty$, $1\leq j, l\leq
N$. Obviously, under Assumptions \ref{asum1}--\ref{asum9}, as
$\mathbf{n}\Rightarrow\infty$, the conclusions of
Theorems \ref{thm2.1} and \ref{thm2.2} hold. Furthermore, in
Corollary~\ref{cor1}, the conditions on $q
=q_{\mathbf{n}}\rightarrow\infty$ and the bandwidth
$h_{\mathbf{n}}$ can be modified as
$\tilde{\mathbf{n}}q^{-\mu}\rightarrow0$,
$q=\mathrm{o}((\tilde{\mathbf{n}}h_{\mathbf{n}}^{k})^{1/2N})$ and
$qh_{\mathbf{n}}^{\delta k/[a(2+\delta)]}>1$ for some
$(4+\delta)N/(2+\delta)<a<\mu\delta/(2+\delta)-N$, then, under
Assumptions \ref{asum1}--\ref{asum7}, the conclusions of Theorems \ref{thm2.1} and
\ref{thm2.2} still hold. Similarly, in Corollary~\ref{cor2},
if $(\tilde{\mathbf{n}}h_{\mathbf{n}}^{k})^{1/2N}
h_{\mathbf{n}}^{\delta
k/[a(2+\delta)]}(\ln\tilde{\mathbf{n}})^{-1} \rightarrow\infty$
for $a>(4+\delta)N/(2+\delta)$, then, under Assumptions \ref{asum1}--\ref{asum7}, the
conclusions of Theorems \ref{thm2.1} and~\ref{thm2.2} still hold.
\end{rem}
%
\begin{rem}\label{rem2}
 If $\chi(n',n'')\leq C_{5}(n'+n''+1)^{\kappa}$
for some $C_{5}>0$ and $\kappa>1$, let the condition
$\tilde{\mathbf{n}}\varphi(q)\rightarrow0$ in Assumption~\ref{asum9} be
replaced by
$(\tilde{\mathbf{n}}^{\kappa+1}/\tilde{\mathbf{p}})\varphi(q)\rightarrow
0$ as $\tilde{\mathbf{n}}\rightarrow\infty$, then the conclusions
of Theorems \ref{thm2.1} and \ref{thm2.2} still hold. In this
case, analogues of Corollaries \ref{cor1} and \ref{cor2} and Remark~\ref{rem1} can also
be obtained.
\end{rem}

\section{Quantile regression and robust smoothers with functional coefficients}\label{sec3}

The general Theorems \ref{thm2.1} and \ref{thm2.2} have different
applications depending on the
choice of $\rho(\cdot)$ function. In this section, we are particularly
discussing the
spatial regression problems with functional coefficients for the
conditional quantiles and robust
functionals.

\subsection{Quantile regression}\label{sec3.1}
Let $F(\cdot|X,U)$ denote the conditional distribution of $Y$
given $X$ and $U$. Then the $\tau$th conditional quantile of $Y$
given $X$ and $U$ is $F^{-1}(\tau|X,U)$, for $0<\tau<1$.
Conditional quantiles have several advantages over conditional
means. For example, they can be defined without any moment
restrictions on $Y$. Plotting the 0.25th, 0.5th, and 0.75th
conditional quantiles would give us more understanding on the data
than plotting just the conditional mean. Quantile regression can
also be useful for the estimation of predictive intervals. For
example, estimates of $F^{-1}(\tau/2|X,U)$ and
$F^{-1}(1-\tau/2|X,U)$ can be used to obtain a $100(1-\tau)\%$
nonparametric interval of prediction of the response given $X$
and $U$. Hallin \textit{et al}. \cite{s15} have studied the spatial
conditional quantile regression estimation, which may however
suffer from curse of dimensionality in general.

We estimate the $\tau$th conditional quantile of $Y_{\mathbf{i}}$
given $X_{\mathbf{i}}$ and $U_{\mathbf{i}}$, approximated by the
functional-coefficient linear structure in (\ref{1.2a}) with
$\beta_r(u)$'s defined by minimizing (\ref{1.2b}) with
$\rho_{\tau}(\cdot)$ instead of $\rho(\cdot)$. If $\tau=1/2$, we
estimate the conditional median. Let
$\hat{a}_{\tau},\hat{b}_{\tau}$ be the minimizer of (\ref{2.1})
with $\rho_{\tau}(\cdot)$ instead of $\rho(\cdot)$. Set
$\hat{\beta}_{\tau}(u_{0})=\hat{a}_{\tau}$. Then the estimator of
the $\tau$th conditional quantile of $Y$ given $X=x$ and
$U=u_{0}$, approximated by the functional-coefficient linear
structure, is $\hat{Y}_{\tau}=x^{T}\hat{\beta}_{\tau}(u_{0})$. To
state the asymptotic results, we need \DIFadd{the following.} %
\renewcommand{\theassum}{Q}
\begin{assum}\label{asumQ}
There exist positive constants
$\bar{c}_{6},C_{6}$ such that the conditional density function
$f_{\varepsilon}(y|X_{i},U_{\mathbf{i}})$ of
$\varepsilon_{\mathbf{i}}$ given $X_{i},U_{\mathbf{i}}$ satisfies
that
$|f_{\varepsilon}(y|X_{i},U_{\mathbf{i}})-f_{\varepsilon
}(0|X_{i},U_{\mathbf{i}})|\leq
C_{6}|y|$ for all $y\in[-\bar{c}_{6},\bar{c}_{6}]$, where
$\varepsilon_{\mathbf{i}}$ is defined in Assumption~\ref{asum3}.
\end{assum}

In this case, since $\psi_{\tau}(z)=2\tau I(z>0)+2(\tau-1)I(z<0)$,
it is easy to show that
Assumption~\ref{asum4} holds and $E(\psi_{\tau}^{2}(\varepsilon_{\mathbf
{i}})|X_{\mathbf{i}},U_{\mathbf{i}})=4\tau(1-\tau)$.
If Assumption~\ref{asumQ} holds, then Assumption~\ref{asum3}
holds with
$\phi(X_{i},U_{\mathbf{i}})=2f_{\varepsilon}(0|X_{i},U_{\mathbf{i}})$.
Let $
\Phi_{\tau}(u)=2E(f_{\varepsilon}(0|X_{i},U_{\mathbf{i}})X_{\mathbf
{i}}X_{\mathbf{i}}^{T}|U_{\mathbf{i}}=u)$
and
$\Omega(u)=E(X_{\mathbf{i}}X_{\mathbf{i}}^{T}|U_{\mathbf{i}}=u)$.
Applying 
Theorems \ref{thm2.1} and \ref{thm2.2} to
quantile regression, we 
have \DIFadd{the following theorem.} 
%
\begin{teo}\label{thm3.1}
\textup{(1)} Assume that Assumptions \ref{asum1}, \ref{asum2}, \ref{asum5}--\ref{asum9}
and
\textup{\ref{asumQ}} hold. Suppose $\Phi_{\tau}(u)$ and
$\Omega(u)$ are continuous in some neighborhood of $u_{0}$ and
$\Phi_{\tau}(u_{0})$ is positive definite, with $u_{0}$ an
interior point of $E_0$. Then, as $\mathbf{n}\rightarrow\infty$,
\[
\sqrt{\tilde{\mathbf{n}}h_{\mathbf{n}}^{k}}\left(\hat{\beta}_{\tau
}(u_{0})-\beta(u_{0})-
\frac{h_{\mathbf{n}}^{2}}{2\mu_{0}}\zeta(u_{0})\right)\rightarrow_{d}
N\biggl(0,
\frac{4\tau(1-\tau)\nu_{0}}{f(u_{0})\mu_{0}^{2}}\Phi_{\tau
}^{-1}(u_{0})\Omega(u_{0})\Phi_{\tau}^{-1}(u_{0})\biggr).
\]

\textup{(2)} Assume that Assumptions \ref{asum1}, \ref{asum2}, \ref{asum5}--\ref{asum9}
and
\textup{\ref{asumQ}} hold in some right neighborhood of
$u_{*}$ and $\Phi_{\tau}(u), \Omega(u)$ are continuous in some
right neighborhood of $u_{*}$. Suppose $\Delta_{c}$ is invertible
and $\Phi_{\tau}(u_{*})$ is positive definite. Then, as
$\mathbf{n}\rightarrow\infty$,
\begin{eqnarray*}
&&\sqrt{\tilde{\mathbf{n}}h_{\mathbf{n}}^{k}}\left(\hat{\beta}_{\tau
}(u_{h})-\beta(u_{h})
-\frac{h_{\mathbf{n}}^{2}}{2}\Biggl[\delta_{11}
\bar{\zeta}(u_{*})+\sum_{l=1}^{k}\delta_{1(l+1)}\bar{\zeta
}^{^{(l)}}(u_{*})\Biggr]\right)
\\
&&\quad \rightarrow_{d} N\left(0, \frac{4\tau(1-\tau)\lambda
_{11}}{f(u_{*})}\Phi_{\tau}^{-1}(u_{*})
\Omega(u_{*})\Phi_{\tau}^{-1}(u_{*})\right).
\end{eqnarray*}
\end{teo}

\subsection{Robust smoothers}\label{sec3.2}
It is known 
that the mean is sensitive to outliers, see \cite{s16}
and \cite{s23}. Since the local average estimator is basically a
mean type estimator, it is also sensitive to outliers. To
robustify this procedure, it is suggested that the function
$\rho(\cdot)$ be chosen so that its first derivative is given by
\[
\psi_{c}(z)=\max\bigl\{-1,\min\{z/c,1\}\bigr\},\qquad c>0,
\]
see \cite{s17}
for interesting discussions. 
We estimate the conditional robust smoother of $Y_{\mathbf{i}}$
given $X_{\mathbf{i}}$ and $U_{\mathbf{i}}$, approximated by
(\ref{1.2a}) with $\beta_r(u)$'s defined by minimizing
(\ref{1.2b}), with the $\rho(z)$ that has the derivative
$\psi_{c}(z)$.
\renewcommand{\theassum}{R}
\begin{assum}\label{asumR}
 The conditional density function
$f_{\varepsilon}(y|X_{i},U_{\mathbf{i}})$ of
$\varepsilon_{\mathbf{i}}$ given $X_{i},U_{\mathbf{i}}$ is
symmetric about 0. There is a positive constant $C_{7}$ such that
$f_{\varepsilon}(y|X_{i},U_{\mathbf{i}})\leq C_{7}$.
\end{assum}

Let $\hat{a}_{c},\hat{b}_{c}$ be the minimizer of (\ref{2.1}) with
the $\rho(\cdot)$ satisfying that
$\rho'(z)=\psi_{c}(z)
,c>0$. Set $\hat{\beta}_{c}(u_{0})=\hat{a}_{c}$. In this case,
Assumption~\ref{asum4} holds automatically. If Assumption~\ref{asumR} holds, then
Assumption~\ref{asum3} holds with
$\phi(X_{i},U_{\mathbf{i}})=P\{|\varepsilon_{\mathbf{i}}|\leq
c|X_{i},U_{\mathbf{i}}\}/c$. Let
$\Phi_{c}(u)=E(P\{|\varepsilon_{\mathbf{i}}|\leq
c|X_{i},U_{\mathbf{i}}\}X_{\mathbf{i}}X_{\mathbf{i}}^{T}|U_{\mathbf{i}}=u)/c$.
An application of Theorems \ref{thm2.1} and \ref{thm2.2} yields
%
\begin{teo}\label{thm3.2}
Assume that Assumptions \ref{asum1}, \ref{asum2}, \ref{asum5}--\ref{asum9}
and
\textup{\ref{asumR}} hold. Suppose $\Phi_{c}(u)$ and
$\Sigma(u)$ are continuous in some neighborhood of $u_{0}$ and
$\Phi_{c}(u_{0})$ is positive definite. Then, as
$\mathbf{n}\rightarrow\infty$, the conclusions of
Theorems \ref{thm2.1} and \ref{thm2.2} hold with $\Phi(u)$
replaced by $\Phi_{c}(u)$.
\end{teo}
%
\begin{rem}\label{rem3}
 Analogues of Theorems \ref{thm3.1} and
\ref{thm3.2} can also be obtained under the conditions that
$\chi(n',n'')\leq C_{5}(n'+n''+1)^{\kappa}$ and (or)
$\mathbf{n}\Rightarrow\infty$ and (or) $\varphi(t)=\mathrm{O}(t^{-\mu})$
for some $\mu>2(3+\delta)N/\delta$ or $\varphi(t)=\mathrm{O}(\mathrm{e}^{-\nu t})$
for some $\nu>0$, details are omitted for the sake of brevity.
\end{rem}

\section{Random fields with a spatial trend}\label{sec4}

In Section~\ref{sec2}, the stationary process
$\{Y_{\mathbf{i}},X_{\mathbf{i}},U_{\mathbf{i}}\}$ was assumed to
be observed. This assumption may often be violated in practice. As
a reasonable alternative, we can assume that nonstationarity is
due to the presence of a spatial trend, as done in \cite{s15}, and that, instead, 
we actually observe
$\{\tilde{Y}_{\mathbf{i}},\tilde{X}_{\mathbf{i}},\tilde{U}_{\mathbf
{i}}\}$,
with
%
\begin{equation}\label{4.1}
\tilde{Y}_{\mathbf{i}}=\alpha_{Y}(s_{\mathbf{i}})+Y_{\mathbf{i}},\qquad
\tilde{X}_{\mathbf{i}}=\alpha_{X}(s_{\mathbf{i}})+X_{\mathbf{i}},\qquad
\tilde{U}_{\mathbf{i}}=\alpha_{U}(s_{\mathbf{i}})+U_{\mathbf{i}},
\end{equation}
where $s_{\mathbf{i}}=(s_{i_{1}},\ldots, s_{i_{N}}) :=
(i_{1}/n_{1},\ldots, i_{N}/n_{N})$ and $s\in[0, 1]^{N}\rightarrow
(\alpha_{Y}(s),\alpha_{X}(s),\alpha_{U}(s))$ is some deterministic
but unknown trend function.

For the sake of simplicity, we assume throughout this section that
$N = 2$, which is the most frequent case in practice. Since
$(Y_{\mathbf{i}},X_{\mathbf{i}},U_{\mathbf{i}})=(\tilde{Y}_{\mathbf
{i}}-\alpha_{Y}(s_{\mathbf{i}}),
\tilde{X}_{\mathbf{i}}-\alpha_{X}(s_{\mathbf{i}}),\tilde{U}_{\mathbf
{i}}-\alpha_{U}(s_{\mathbf{i}}))$
is 
unobservable, 
the analysis proceeds in two steps. First, obtain an estimation of
the spatial trend $(\alpha_{Y}(s_{\mathbf{i}}),
\alpha_{X}(s_{\mathbf{i}}),\alpha_{U}(s_{\mathbf{i}}))$ via
kernel smoothing method. In the second step, the detrended data is
supposed to satisfy the stationarity assumption, 
yielding the estimated coefficient function $\check{\beta}_{r}(u),
r=1,\ldots,d$ with the detrended ${Y}_{\mathbf{i}}$'s,
${X}_{\mathbf{i}}$'s and ${U}_{\mathbf{i}}$'s.

Let
\[
w(s_{\mathbf{i}},s)=\frac{W((s_{\mathbf{i}}-s)/g_{\mathbf{n}})}{\sum
_{\mathbf{j}\in
G_{N}}W((s_{\mathbf{j}}-s)/g_{\mathbf{n}})},
\]
where $g_{\mathbf{n}}$ is a bandwidth tending to $0$ and
$W(\cdot)$ is a chosen kernel function. Then the kernel estimators
of $\alpha_{Y}(s), \alpha_{X}(s)$ and $\alpha_{U}(s)$ are
%
\begin{equation}\label{4.2}
\hat{\alpha}_{Y}(s)=\sum_{\mathbf{i}\in
G_{N}}\tilde{Y}_{\mathbf{i}}w(s_{\mathbf{i}},s), \qquad \hat{\alpha
}_{X}(s)=\sum_{\mathbf{i}\in
G_{N}}\tilde{X}_{\mathbf{i}}w(s_{\mathbf{i}},s), \qquad \hat{\alpha
}_{U}(s)=\sum_{\mathbf{i}\in
G_{N}}\tilde{U}_{\mathbf{i}}w(s_{\mathbf{i}},s).
\end{equation}
Let
$\hat{Y}_{\mathbf{i}}=\tilde{Y}_{\mathbf{i}}-\hat{\alpha}_{Y}(s_{\mathbf{i}})$,
$\hat{X}_{\mathbf{i}}=\tilde{X}_{\mathbf{i}}-\hat{\alpha}_{X}(s_{\mathbf{i}})$
and
$\hat{U}_{\mathbf{i}}=\tilde{U}_{\mathbf{i}}-\hat{\alpha}_{U}(s_{\mathbf{i}})$.
Based on the estimated spatial data
$\{(\hat{Y}_{\mathbf{i}},\hat{X}_{\mathbf{i}},\hat{U}_{\mathbf{i}})\dvtx
\mathbf{i}\in
G_{N}\} $,
we solve the following minimization problem
%
\begin{equation}\label{eq4.3}
\min_{a_{r},b_{r},r=1,\ldots,d}\sum_{\mathbf{i}\in
G_{N}}\rho\left(\hat{Y}_{\mathbf{i}}-\sum_{r=1}^{d}\bigl[a_{r}+(\hat
{U}_{\mathbf{i}}-u_{0})^{T}b_{r}\bigr]\hat{X}_{\mathbf{i}r}\right)
K\biggl(\frac{\hat{U}_{\mathbf{i}}-u_{0}}{h_{\mathbf{n}}}\biggr).
\end{equation}
Let
$\check{a}_{r},\check{b}_{r}$ be the minimizer of (\ref{eq4.3}). Set
$\check{a}=(\check{a}_{1},\ldots,\check{a}_{d})^{T}$. Then the
M-estimator of
$\beta(u_{0})=(\beta_{1}(u_{0}),\ldots,\beta_{d}(u_{0}))^{T}$ is
\[
\check{\beta}(u_{0})
=\bigl(\check{\beta}_{1}(u_{0}),\ldots,\check{\beta}_{d}(u_{0})\bigr)^{T}=\check{a}.
\]

To study the asymptotic behavior of the new estimators, we need
the following additional 
conditions similar to those in \cite{s15}.
\begin{enumerate}[(B5)]
\item[(B1)] $E|Y_{\mathbf{i}}|^{2+\delta}<\infty$,
$E\|X_{\mathbf{i}}\|^{2+\delta}<\infty$ and
$E\|U_{\mathbf{i}}\|^{2+\delta}<\infty$ for some $\delta>0$ and
$\varphi(z)$ in (\ref{1.1})
satisfies that $\varphi(z)<\tilde{C}_{0}z^{-\beta}$, where
$0<\tilde{C}_{0}<\infty$ and $\beta>(1+(1+\delta)(1+N))/\delta$.

\item[(B2)] For
$\varrho=(\beta-1-N-(1+\beta)/(1+\delta))/(\beta+3-N-(1+\beta)/(1+\delta))$,
$\ln\tilde{\mathbf{n}}/(\tilde{\mathbf{n}}^{\varrho}g_{\mathbf{n}}^{N})=\mathrm{o}(1)$.

\item[(B3)] $s\rightarrow\alpha_{Y}(s)$, $s\rightarrow\alpha_{X}(s)$
and $s\rightarrow\alpha_{U}(s)$ are $m$ times differentiable with
bounded derivatives on $S := [0, 1]^{2}$, where $m$ is some
positive integer.

\item[(B4)] There exists a continuous sampling intensity (density)
function $\tilde{f}$ defined on $S$ and constants $\tilde{c}_{0}$
and $\tilde{c}_{1}$ such that $0 < \tilde{c}_{0}\leq
\tilde{f}(s)\leq\tilde{c}_{1}<\infty$ for any $s\in S$ and
$\tilde{\mathbf{n}}^{-1}\sum_{\mathbf{i}\in
G_{N}}I(s_{\mathbf{i}}\in A)\rightarrow\int_{A}\tilde{f}(s)\,\mathrm{d}s$
for any measurable set $A\subset S$, as
$\tilde{\mathbf{n}}\rightarrow\infty$.

\item[(B5)] The kernel $W(s)$, defined on $R^{2}$, has bounded support
with Lipschitz property, that is $|W(s)-W(s')|\leq
\tilde{C}_{1}\|s-s'\|$ for all $s, s'\in R^{2}$, where
$\tilde{C}_{0}$ is a generic positive constant, and satisfies ($s^{\otimes i}$ stands for the $i$th Kronecker power of $s$)
\begin{eqnarray*}
\int W(s)\,\mathrm{d}s&=&1, \qquad \int s^{\otimes i}W(s)\,\mathrm{d}s=0,\qquad  i=1, \ldots,
m-1,\\
\int s^{\otimes m}W(s)\,\mathrm{d}s&\neq&0.
\end{eqnarray*}
\end{enumerate}

Assumptions (B1) and (B2) are technical conditions for deriving
the convergence of this kernel smoothing; see \cite{s18} for
similar assumptions. Assumption (B4) is mentioned for the sake of
generality, and is trivially satisfied in the case of a regular
grid. Assumptions (B3) and (B5) are standard assumptions on the
smoothness of spatial trend functions and a higher order kernel
function, respectively, which ensure that the bias term of the
spatial trend estimators is of order $\mathrm{O}(g_{\mathbf{n}}^{m})$
(which can also be achieved by a local polynomial fitting of order
$(m-1)$).

We further need to strengthen Assumptions \ref{asum4}--\ref{asum7} as the follows.
\renewcommand{\theassum}{\arabic{assum}$^\prime$}
\setcounter{assum}{3}
\begin{assum}\label{asum4'}
Let $\mathcal{L}_{p}(\mathcal{F})$ denote the
class of $\mathcal{F}$-measurable random variable $\xi$ satisfying
$\|\xi\|_{p}=(E|\xi|^{p})^{1/p}<\infty$. The function
$\psi(\cdot)$ satisfies that
$E(|\psi(\eta_{\mathbf{i}}+\xi)-\psi(\eta_{\mathbf{i}})||X_{\mathbf
{i}},U_{\mathbf{i}})\leq
\tilde{C}_{1}\epsilon$ for 
$\eta_{\mathbf{i}}\in\mathcal{L}_{1}(\mathcal{B}(\{\mathbf{i}\}))$
and $\xi\in\mathcal{L}_{1}(\mathcal{B}(G_{N}))$ such that
$|\xi|<\epsilon$, and that $|\psi(v+s)-\psi(v)|\leq
\tilde{C}_{1}$ for any $|s|\leq c_{2}$ and $v\in R^{1}$, where
$\tilde{C}_{1}$, $\epsilon$ and $c_2$ are some positive constants.
\end{assum}
%
\begin{assum}\label{asum5'}
The bandwidths $h_{\mathbf{n}}$ and
$g_{\mathbf{n}}$ satisfy that $h_{\mathbf{n}}\leq
\tilde{C}_{2}\tilde{\mathbf{n}}^{-1/(k+4)}$ for some positive
constant $\tilde{C}_{2}$,
$g_{\mathbf{n}}^{m}/h_{\mathbf{n}}\rightarrow0$,
$\tilde{\mathbf{n}}h_{\mathbf{\mathbf{n}}}^{k}g_{\mathbf
{n}}^{2m}\rightarrow
0$,
$h_{\mathbf{n}}^{k}\ln\tilde{\mathbf{n}}/g_{\mathbf{n}}^{2}\rightarrow
0$ and
$\ln\tilde{\mathbf{n}}/(\tilde{\mathbf{n}}g_{\mathbf{n}}^{2}h_{\mathbf
{n}}^{2})\rightarrow
0$.
\end{assum}
%
\begin{assum}\label{asum6'}
 $\max_{\mathbf{i}\in
G_{N}}\|X_{\mathbf{i}}\|=\mathrm{O}_{p}(1)$ and
$E\|X_{\mathbf{i}}\|^{4+2\delta}
<\infty$.
\end{assum}
%
\begin{assum}\label{asum7'}
 The kernel function $K(\cdot)\geq0$
is a bounded symmetric function with a compact support $\tilde{M}$
and is continuously differentiable in $\tilde{M}=(-M_{1},M_{1})\times
\cdots\times(-M_{k},M_{k})$
and $\int_{\tilde{M}}uu^{T}K(u)\,\mathrm{d}u$ is positive definite.
\end{assum}

Assumption~\ref{asum4'} is easily checked. 
For example, it holds when $\rho(z)=\rho_{\tau}(z)$ or $\psi(z)=\psi
_c(z)$ and the conditional density of
$\eta_{\mathbf{i}}$ given $X_{\mathbf{i}}$ and $U_{\mathbf{i}}$ is
bounded on $[-\epsilon,\epsilon]$. 
Assumption~\ref{asum5'} on the
bandwidths $h_{\mathbf{n}}$ and $g_{\mathbf{n}}$ is easily
satisfied, and can be weaken as:
$\tilde{\mathbf{n}}^{2/(k+4)}g_{\mathbf{n}}^{m}\rightarrow0$ and
$\ln\tilde{\mathbf{n}}/(\tilde{\mathbf{n}}^{k/(k+4)}g_{\mathbf
{n}}^{2})\rightarrow
0$ if we take the optimal
$h_{\mathbf{n}}=h_{0}\tilde{\mathbf{n}}^{-1/(k+4)}$ for some
$h_{0}>0$. The condition, $\max_{\mathbf{i}\in
G_{N}}\|X_{\mathbf{i}}\|=\mathrm{O}_{p}(1)$, in Assumption~\ref{asum6'} is only a
technical condition, and can also be weakened with
$h_{\mathbf{n}}$ and $g_{\mathbf{n}}$ properly chosen.

We state the asymptotic distribution of the
estimators $\check{\beta}_{r}(u_{0}), r=1,\ldots,d$, as follows.
%
\begin{teo}\label{thm4.1}
Assume that Assumptions \ref{asum1}--\ref{asum3}, \ref{asum4'}--\ref{asum7'}, \ref{asum8}--\ref{asum9}
and \textup{(B1)}--\textup{(B5)} hold and $\Phi(u), \Sigma(u)$ are continuous in some
neighborhood of $u_{0}$ and $\Phi(u_{0})$ is positive definite. If
$u_{0}$ is an interior point of the support of the design density
$f(u)$, then, as $\mathbf{n}\rightarrow\infty$,
\[
\sqrt{\tilde{\mathbf{n}}h_{\mathbf{n}}^{k}}\left(\check{\beta
}(u_{0})-\beta(u_{0})
-\frac{h_{\mathbf{n}}^{2}}{2\mu_{0}}\zeta(u_{0})\right)\rightarrow_{d}
N\biggl(0,
\frac{\nu_{0}}{f(u_{0})\mu_{0}^{2}}\Phi^{-1}(u_{0})\Sigma(u_{0})\Phi
^{-1}(u_{0})\biggr).
\]
\end{teo}

With $\rho_{\tau}(z)$ instead of $\rho(z)$ and Assumption~\ref{asum4'}
replaced by Assumption~\ref{asumQ}, we have
\DIFadd{the following theorem.}
%
\begin{teo}\label{thm4.2}
Assume that Assumptions \ref{asum1}--\ref{asum2}, \ref{asum5'}--\ref{asum7'}, \ref{asum8}--\ref{asum9},
\textup{(B1)}--\textup{(B5)} and
\textup{\ref{asumQ}} hold. Suppose $\Phi_{\tau}(u)$ and $\Omega(u)$ are
continuous in some neighborhood of $u_{0}$ and
$\Phi_{\tau}(u_{0})$ is positive definite and
$f_{\varepsilon}(0|X_{\mathbf{i}},U_{\mathbf{i}})\leq C$ for some
$C>0$. If $u_{0}$ is an interior point of the support of the
design density $f(u)$, then, as $\mathbf{n}\rightarrow\infty$,
\[
\sqrt{\tilde{\mathbf{n}}h_{\mathbf{n}}^{k}}\left(\check{\beta}_{\tau}(u_{0})-
\beta(u_{0})-\frac{h_{\mathbf{n}}^{2}}{2\mu_{0}}\zeta(u_{0})\right
)\rightarrow_{d}
N\biggl(0,
\frac{4\tau(1-\tau)\nu_{0}}{f(u_{0})\mu_{0}^{2}}\Phi_{\tau
}^{-1}(u_{0})\Omega(u_{0})\Phi_{\tau}^{-1}(u_{0})\biggr).
\]
\end{teo}

Similarly, with $\psi_{c}(z)=\max\{-1,\min\{z/c,1\}\},c>0$, it is
easy to check that Assumption~\ref{asum4'} holds. In this case, we have
\DIFadd{the following.}
%
\begin{teo}\label{thm4.3} Assume that Assumptions \ref{asum1}--\ref{asum2}, \ref{asum5'}--\ref{asum7'}, \ref{asum8}--\ref{asum9},
\textup{(B1)}--\textup{(B5)} and
\textup{\ref{asumR}} hold. Suppose $\Phi_{c}(u)$ and
$\Sigma(u)$ are continuous in some neighborhood of $u_{0}$ and
$\Phi_{c}(u_{0})$ is positive definite. Then 
the conclusions of Theorem~\ref{thm4.1} hold with $\Phi(u)$
replaced by $\Phi_{c}(u)$.
\end{teo}

The proofs of Theorems \ref{thm4.1}--\ref{thm4.3} are postponed to
\hyperref[appendix]{Appendix}.

\section{An application to soil data analysis}\label{sec5}

We are analysing a spatial soil data set, soil250, in R package
GeoR, which consists in uniformity trial with 250 undisturbed soil
samples collected at 25 cm soil depth of spacing of 5 meters,
resulting on a regular grid of $25 \times10$ points. The data
frame is with 250 observations on the 22 variables about several
soil chemistry properties measured on the grid. In this analysis,
for simplicity, we only consider 10 variables, which are Linha
(x-coordinate), Coluna (y-coordinate), pHKCl (soil pH by KCl), Ca
(calcium content), Mg (magnesium content), K (potassio content),
Al (aluminium content), C (carbon content), N (nitrogen content),
and CTC (catium exchange capability). 
\begin{figure}

\includegraphics{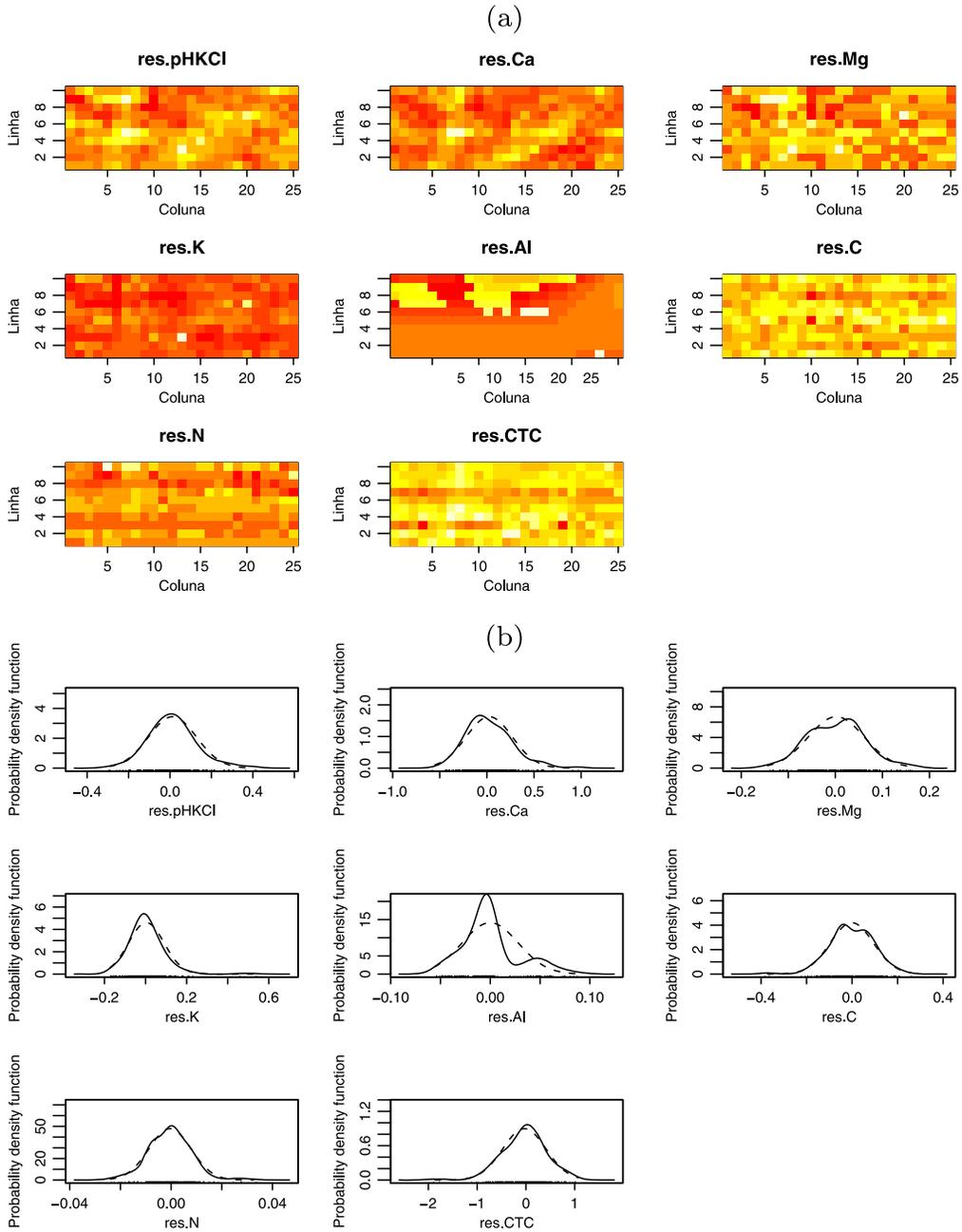}

\caption{Soil data: (a) The images of 8 soil properties variables
after spatial trend removal by sm.regression, plotted over space
(Linha, Coluna);
(b) The kernel density estimates (solid line) of the 8 soil properties variables
after spatial trend removal by sm.regression, where the dashed line is
for the
Gaussian density with the same mean and variance, respectively, for
each variable.}\label{soil2}
\end{figure}
Zheng \textit{et al}. \cite{s45} recently analysed the spatial spectral
density for the CTC variable. Our objective here is to analyse
the impacts of the soil chemistry properties of Ca, Mg, K, Al, C
and N as well as the soil chemistry property index pHKCl on the
CTC, an important soil property for soil reservation concerned
with in agriculture
science. 
In the original data, there seem to be some spatial trends for all
variables, so we 
apply sm.regression in R package sm to remove the spatial trends.
The resulting spatial data of these variables, denoted by prefix
``res.'' standing for residual, are plotted in panel (a) of
Figure~\ref{soil2}, which appear more stationary. We hence analyse
the
relationship of these variables based on the residual data, 
the kernel density estimates (in solid line) of which are also plotted
in Figure~\ref{soil2}(b), with the
dashed line for the Gaussian density of the same mean and
variance. It is clear that the distribution of the response,
res.CTC, is quite close to normal, indicating that mean and
median regression analyses are similar (only median regression is
provided below). Further, considering the spatial neighbouring
effects, we also include the nearest neighbour variables of the
CTC, denoted by res.CTCw, res.CTCe, res.CTCn and res.CTCs for the
west, east, north and south nearest neighbours. Thus we have 11
covariates, including res.pHKCl and the soil chemistry property
variables (res.Ca, res.Mg, res.K, res.Al, res.C and res.N) as well
as the four neighbouring variables of the CTC. It is impossible to
apply general nonparametric quantile analysis of the impacts of
these covariates on the response as done in \cite{s15} as it suffers from severe ``curse of dimensionality''.

To have a preliminary understanding of the possible relationship,
we made a simple nonparametric regression analysis of the CTC
on each covariate by applying sm.regression in R package sm
(the results not reported here to save space).
It
appears that
the response res.CTC is basically linearly related with each of the individual
covariates, suggesting we may 
consider regressing $Y_{ij}=\mathrm{res.CTC}_{ij}$ at a grid $(i,j)$
on the covariates in a linear
form
%
\begin{eqnarray}\label{soil-FC}
&&a_0(\mathrm{res.pHKCl}_{ij})+a_1(\mathrm{res.pHKCl}_{ij}) \mathrm
{res.Ca}_{ij}+a_2(\mathrm{res.pHKCl}_{ij}) \mathrm{res.Mg}_{ij}
\nonumber\\
&&\quad {}+a_3(\mathrm{res.pHKCl}_{ij}) \mathrm{res.K}_{ij}
+a_4(\mathrm{res.pHKCl}_{ij}) \mathrm{res.Al}_{ij}+a_5(\mathrm
{res.pHKCl}_{ij}) \mathrm{res.C}_{ij} \nonumber\qquad\\
&&\quad {}+a_6(\mathrm{res.pHKCl}_{ij}) \mathrm{res.N}_{ij}+a_7(\mathrm
{res.pHKCl}_{ij}) \mathrm{res.CTCw}_{ij} \\
&&\quad {} +a_8(\mathrm{res.pHKCl}_{ij})
\mathrm{res.CTCe}_{ij}+a_9(\mathrm{res.pHKCl}_{ij}) \mathrm{res.CTCn}_{ij}\nonumber\\
&&\quad {} +a_{10}(\mathrm
{res.pHKCl}_{ij}) \mathrm{res.CTCs}_{ij}\nonumber
\end{eqnarray}
for $ 1\le i\le25, 1\le j\le10$, where we take the chemical
property index variable, $\mathrm{res.pHKCl}$, as a regime
variable $U$
and are
interested in the effects of this
index variable in the coefficient functions $a_1(\cdot), \ldots,
a_{10}(\cdot)$ of the components of $X$ denoted for the vector of
other variables, for example, whether these coefficient functions
are constant or not. Here $a_0(\cdot)$ is the baseline effect from
the index variable.

\begin{figure}

\includegraphics{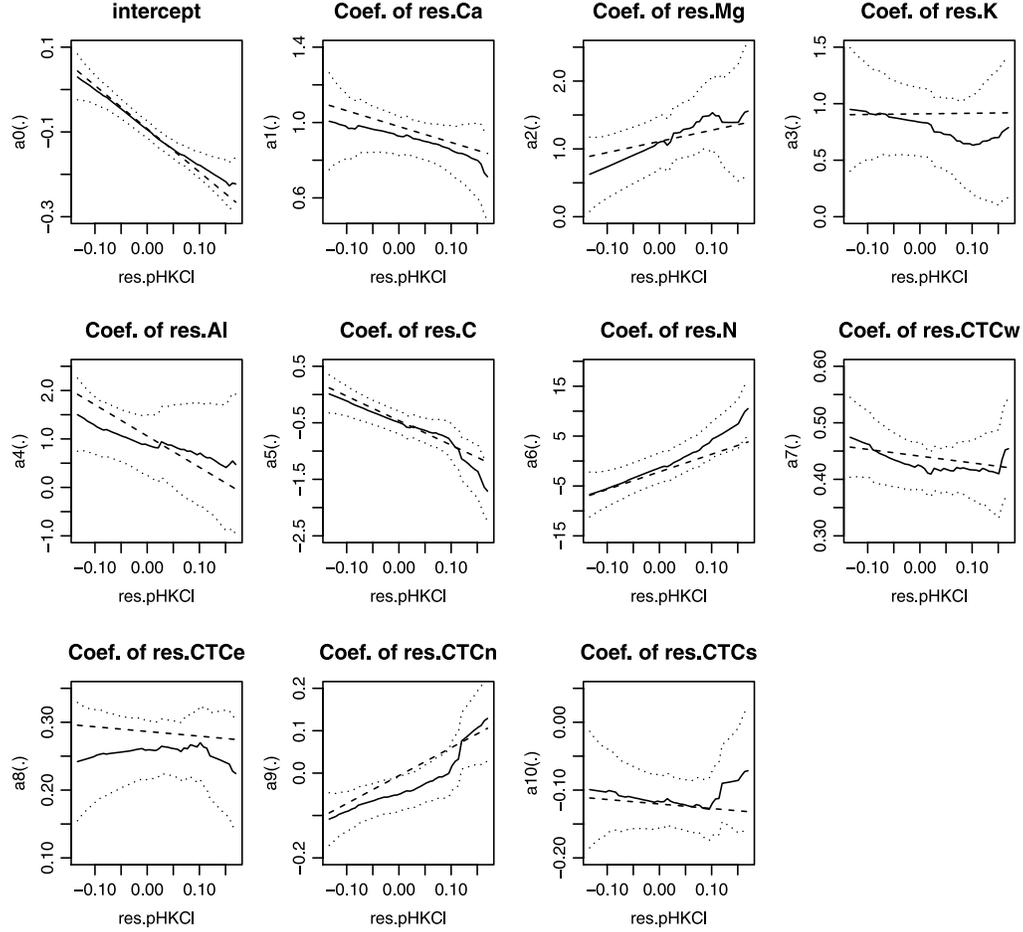}

\caption{Soil data: The median regression estimate ($\tau=0.5$) of the
functional coefficients in (\protect\ref{soil-FC}). The
solid (---) line is for the quantile regression with nonparametric
functional coefficients in this paper by using the selected bandwidth
of 0.263 
by the leave-one-out 
CV, the
dashed ($--$) line is for the
functional coefficient of parametric linear function, and the dotted
($\cdots$) lines are for the 95\% confidence intervals
constructed by asymptotic normality.}\label{soil4}
\end{figure}

\begin{figure}

\includegraphics{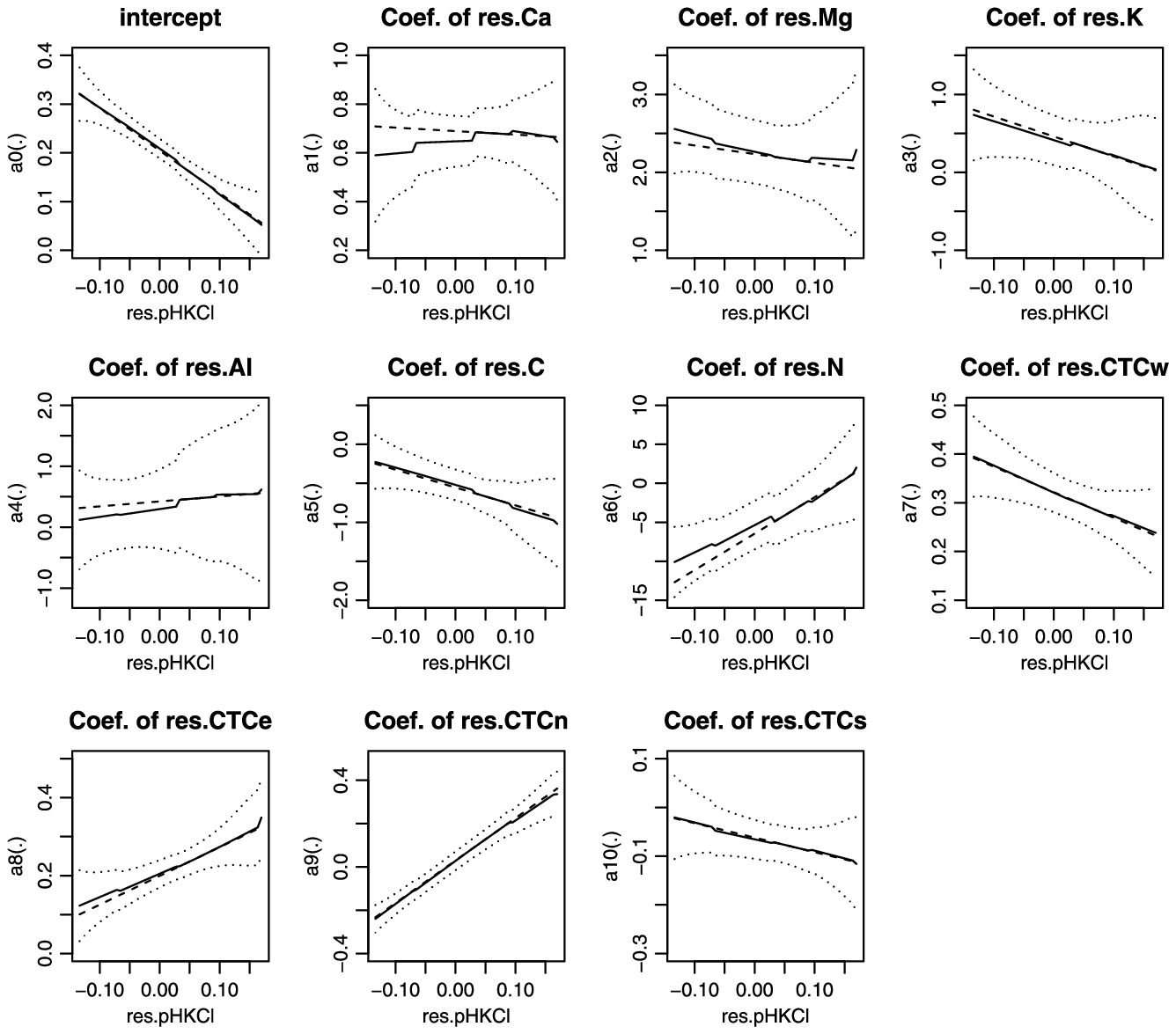}

\caption{Soil data: The quantile regression estimate ($\tau=0.85$) of
the functional coefficients in (\protect\ref{soil-FC}). The
solid (---) line is for the quantile regression with nonparametric
functional coefficients in this paper by using the selected bandwidth
of 0.263 
by the leave-one-out 
CV, the
dashed ($--$) line is for the
functional coefficient of parametric linear function, and the dotted
($\cdots$) lines are for the 95\% confidence
intervals constructed by asymptotic normality.}\label{soil5}
\end{figure}

\begin{figure}

\includegraphics{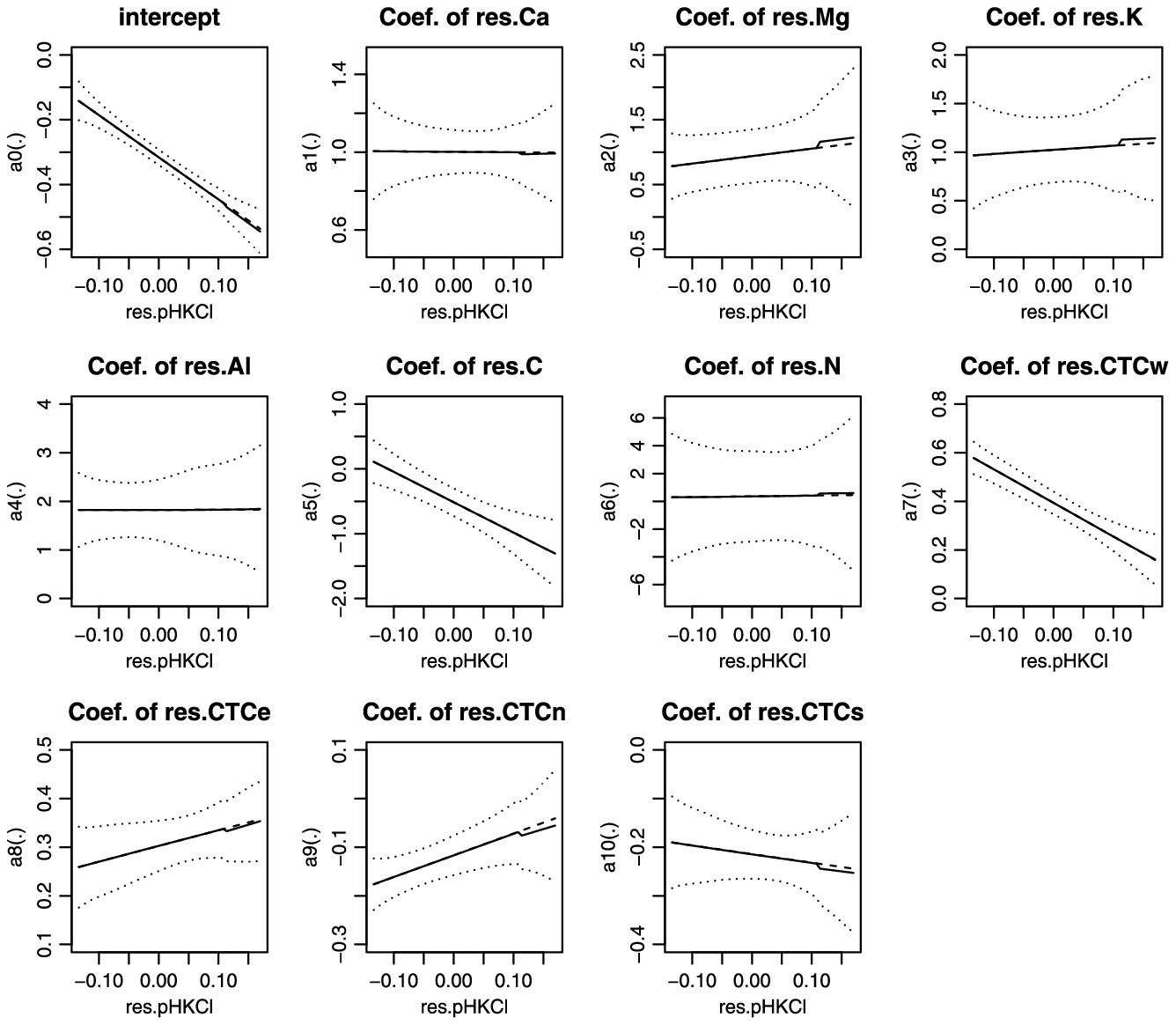}

\caption{Soil data: The quantile regression estimate ($\tau=0.15$) of
the functional coefficients in (\protect\ref{soil-FC}). The
solid (---) line is for the quantile regression with nonparametric
functional coefficients in this paper by using the selected bandwidth
of 0.263 
by the leave-one-out 
CV, the
dashed ($--$) line is for the
functional coefficient of parametric linear function, and the dotted
($\cdots $) lines are for the 95\% confidence intervals constructed by
asymptotic normality.}\label{soil6}
\end{figure}

We here suggest selecting the required bandwidth $h$ in
(\ref{2.1}) by applying an empirical rule of Fan \textit{et al}.
\cite{s7} with cross-validation (CV) of Stone \cite{66} using the check
function $\rho(z)$ in (\ref{2.1}). In the time series context, the
argument for cross-validation as an appropriate method for the
bandwidth can be found in \cite{s24,35a} and \cite{s43a}, among others. This empirical rule of
bandwidth selection procedure is computationally efficient \cite{s7}; see also \cite{s30}, \DIFadd
{Section} 2.3,
in the least squares setting. We first examine the median
regression under $\tau=0.5$, with the range of $h$ taken between
0.15 and 0.3 (partitioned into $q$ small intervals of length
$0.01$). The spatial quantile estimates of these coefficient
functions under $\tau=0.5$ are provided in Figure~\ref{soil4} in
solid lines, for the selected bandwidth of $h=0.263$ by a
leave-one-out CV with $\rho(z)=|z|$ in (\ref{2.1}). In order to
take into account the dependence in the observations, we also
applied a leave-five-out CV for the selection of bandwidth with
$h=0.235$ selected, by which the estimated median regression
coefficient functions are very similar to those with
leave-one-out CV, and are therefore omitted in Figure~\ref{soil4}.
It seems that many of the functional coefficients, such as the
baseline function $a_0(\cdot)$, are nearly linear. We hence also
made median
regression analysis with the coefficient functions 
of 
a linear form, reported in dashed lines in
Figure~\ref{soil4}. In order to examine the impacts of the covariates
on the high or
low CTC variable, we also made similar analysis of spatial
quantile regression of (\ref{soil-FC}) under $\tau=0.85$ and
$\tau=0.15$,
plotted in Figure~\ref{soil5} and Figure~\ref{soil6}, respectively.
In view of the sparsity of extreme data, the range of $h$ was
taken a bit larger between $0.25$ and $0.6$ (with refined
partition of $q$ small intervals of length $0.001$), with the
leave-one-out
CV-selected bandwidths equal to 0.5 
and 0.487 
for $\tau=0.85$ and
$\tau=0.15$, respectively. Again the estimated coefficient
functions based on leave-five-out CV, which are omitted here, are
similar to those with leave-one-out CV.

As the information on how
variable the estimates are would be interesting for statistical
inference, we have also provided pointwise bands, that is\DIFadd{,} a
collection of confidence intervals, for the quantile coefficient
estimates on the basis of the asymptotic theorem
(Theorem~\ref{thm3.1}), which are plotted in dotted lines in
Figures \ref{soil4}--\ref{soil6}. Here the key difficulty in doing
so is the estimation of $f_{\varepsilon}(0|X_{i},U_{\mathbf{i}})$
associated with $\Phi_{\tau}(u)$ in the asymptotic variance of
Theorem~\ref{thm3.1}. Note that we cannot simply assume
$\varepsilon_{\mathbf{i}}$ and $(X_{\mathbf{i}},U_{\mathbf{i}})$
are independent as in the traditional varying-coefficient analysis
in the literature. Therefore\DIFadd{,} the estimation suffers from severe
curse of dimensionality (note that the dimension of
$(X_{\mathbf{i}},U_{\mathbf{i}})$ is equal to 11) at a first
glance. Fortunately, however, by applying indepTest in the R
package ``copula'' with the independence test method of Genest and R\'{e}millard \cite{s12gr}, we find at $5\%$ significance level that the
estimated $\varepsilon_{\mathbf{i}}$ is only dependent on
$(X_{\mathbf{i},8},X_{\mathbf{i},9})=(\mathrm{res.CTCw}, \mathrm{res.CTCe})$ at
$\tau=0.85$, and on
$(X_{\mathbf{i},10},X_{\mathbf{i},11})=(\mathrm{res.CTCn}, \mathrm{res.CTCs})$ at
$\tau=0.15$, while the estimated $\varepsilon_{\mathbf{i}}$ and
$(X_{\mathbf{i}},U_{\mathbf{i}})$ are independent at $\tau=0.5$.
Hence\DIFadd{,} we can easily estimate the conditional density function
$f_{\varepsilon}(0|X_{i},U_{\mathbf{i}})$ by applying npcdens in
the R package ``np'' with the method of Hall \textit{et al}. \cite{s13}.
The asymptotic variance of Theorem~\ref{thm3.1} can thus be
calculated, by which the confidence intervals are constructed.

Let us first examine Figure~\ref{soil6} with $\tau=0.15$. We
can see in this figure that the coefficient functions are close to linear
lines, which, except $a_6(\cdot)$,
are significant from zero at $5\%$ significance level.
For Figure~\ref{soil5} with $\tau=0.85$, though the coefficient
functions are also close to be linear, the
magnitudes of the effects of different covariates through the regime index,
res.pHKCl, appear quite significantly different from those in
Figure~\ref{soil6}.
These findings are also different from that in
Figure~\ref{soil4} with $\tau=0.5$,
meaning that the
effects of the different covariates through the regime index,
res.pHKCl, on CTC perform differently at low, median and high
values of CTC:
\begin{enumerate}[(3)]
\item[(1)] Different covariate effects: The covariate effects 
under different $\tau$'s appear quite different in magnitude, but
mostly are of the same signs. Here res.Ca, res.Mg, res.K and res.Al
have nonnegative effects
for which we cannot reject their constancy, while
res.C has a negative effect decreasing with res.pHKCl, at $5\%$ level
of significance.
However, for the
covariate res.N, it seems clear at $5\%$ significance level that
$a_6(\cdot)$ is not significant from zero under $\tau=0.15$, but
it is an increasing function that is negative (turning to positive
values) when the regime, res.pHKCl,
is less than the thresholds, 0.05 and 0.10, under
$\tau=0.5$ and $\tau=0.85$, respectively.
It looks that the chemistry properties of N (nitrogen content)
may play a
significantly different
role with the regime in the adjustment of the high/low CTC in the
soil. These findings are beyond the traditional median
or mean regression analysis.

\item[(2)] Different neighbouring effects: The neighbouring effects in
quantile analysis
under different $\tau$'s also appear different in magnitude,
and
mostly are of the same sign in the coefficients (positively correlated
with west and east neighbours but negatively with south).
However, it looks at $5\%$ significance level that the CTC in the soil
has a negative correlation with its north neighbour res.CTCn
(note the coefficient $a_9(\cdot)$ is negative) under $\tau=0.15$,
but becomes positively correlated with its north neighbour when the
regime, res.pHKCl,
is over the thresholds 0.10 and 0.05 under $\tau=0.5$ and $\tau=0.85$,
respectively.

\item[(3)] Different regime effects: The regime effects of res.pHKCl seem more
involved under $\tau=0.5$, in particular in the coefficients
of res.C, res.N and res.CTCn, which appear marginally nonlinear at $5\%
$ significance level. For high ($\tau=0.85$) and
low ($\tau=0.15$)
quantiles, the regime effects of res.pHKCl appear linear or constant.
\end{enumerate}

To sum up, although the above analysis is illustrative only, it seems
apparent that the functional-coefficient
spatial quantile regression proposed in this paper is helpful to
uncover and understand the underlying relationship of the soil
chemistry properties with CTC (catium exchange capability) through
the regime index pHKCl. These properties are interesting and
important topics in soil reservation and management.

\begin{appendix}\label{appendix}
\section*{Appendix: Proofs}\label{Aproofs}

In this section, we only sketch the proof of the main theorems
with the necessary lemmas listed. The detail of the proof of the
lemmas is 
much more
complicated and we describe it in detail in the supplementary material
 \cite{s31a}.

Let $C$ denote a generic positive constants not depending on
$\mathbf{n}$, which may take on different values at each
appearance. Under Assumption~\ref{asum2}, by Taylor expansion, for
$U_{\mathbf{i}}=(U_{\mathbf{i}1},\ldots,U_{\mathbf{i}k})^{T}$ such
that
$|U_{\mathbf{i}l}-u_{0l}|\leq Mh_{\mathbf{n}}, 1\leq l\leq k$, we have
$\beta_{r}(U_{\mathbf{i}})=\beta_{r}(u_{0})+\dot{\beta
}_{r}(u_{0})^{T}(U_{\mathbf{i}}
-u_{0})+\frac{1}{2}(U_{\mathbf{i}}-u_{0})^{T}\ddot{\beta}_{r}(\xi
_{\mathbf{i}r})(U_{\mathbf{i}}-u_{0})$,
where
$\dot{\beta}_{r}(u_{0})=(\dot{\beta}_{r1}(u_{0}),\ldots,\dot{\beta
}_{rk}(u_{0}))^{T}$
stands for the gradient of $\beta_{r}(u)$ with respect to $u$ at
$u=u_0$, and
$\xi_{\mathbf{i}r}=(\xi_{\mathbf{i}r1},\ldots,\xi_{\mathbf{i}rk})^{T}$
satisfies that
$|\xi_{\mathbf{i}rl}-u_{0l}|<|U_{\mathbf{i}l}-u_{0l}|$ for $1\leq
l\leq k$.
Denote $\beta^{\ast}(\xi_{\mathbf{i}})=((U_{\mathbf{i}}-u_{0})^{T}\ddot
{\beta}_{1}(\xi_{\mathbf{i}1})(U_{\mathbf{i}}
-u_{0}),\ldots,(U_{\mathbf{i}}-u_{0})^{T}\ddot{\beta}_{d}(\xi_{\mathbf
{i}d})(U_{\mathbf{i}}-u_{0}))^{T}$,
and $e_{\mathbf{i}}
=\frac{1}{2}\beta^{\ast}(\xi_{\mathbf{i}})^{T}X_{\mathbf{i}}$.
Let $\tilde{b}_{l}=(b_{1l},\ldots,b_{dl})^{T}$,
$\tilde{\beta}_{l}(u_{0})=(\dot{\beta}_{1l}(u_{0}),\ldots,\dot{\beta
}_{dl}(u_{0}))^{T}$,
$\tilde{a}=(a_{1},\ldots,a_{d})^{T}$, $\tilde{b}=(\tilde
{b}_{1}^{T},\ldots,\tilde{b}_{k}^{T})^{T}$,
$\dot{\beta}(u_{0})=(\tilde{\beta}_{1}(u_{0})^{T},\ldots
,\break \tilde{\beta}_{k}(u_{0})^{T})^{T}$.
Set $Z_{\mathbf{i}}=(\widetilde{\mathbf{n}}h_{\mathbf
{n}}^{k})^{-1/2}(1,h_{\mathbf{n}}^{-1}(U_{\mathbf{i}}
- u_{0})^{T})^{T}\otimes X_{\mathbf{i}}$ and
$t=(\widetilde{\mathbf{n}}h_{\mathbf{n}}^{k})^{1/2}((\tilde{a}-\beta
(u_{0}))^{T},\break h_{\mathbf{n}}(\tilde{b}
-\dot{\beta}(u_{0}))^{T})^{T}$, where $\otimes$ is the Kronecker
product. Then we have the following new optimization problem
\setcounter{equation}{0}
\begin{eqnarray}\label{A.1}
\hat{t}=\mathop{\Argmin}_{t}\sum_{\mathbf{i}\in
G_{N}}\bigl[\rho\bigl(\varepsilon_{\mathbf{i}}+e_{\mathbf{i}}-t^{T}Z_{\mathbf{i}}\bigr)-
\rho(\varepsilon_{\mathbf{i}}+e_{\mathbf{i}})\bigr]K\biggl(\frac{U_{\mathbf
{i}}-u_{0}}{h_{\mathbf{n}}}\biggr).
\end{eqnarray}
Clearly
%
\begin{eqnarray}\label{A.2}
\hat{t}=\bigl(\widetilde{\mathbf{n}}h_{\mathbf{n}}^{k}\bigr)^{1/2}\bigl(\bigl(\hat{a}-\beta
(u_{0})\bigr)^{T},
h_{\mathbf{n}}\bigl(\hat{b}-\dot{\beta}(u_{0})\bigr)^{T}\bigr)^{T}.
\end{eqnarray}

Denote the objective function in (\ref{A.1}) by $S_{\mathbf{n}}(t)$
and set
\[
\Gamma_{\mathbf{n}}(t)=\sum_{\mathbf{i}\in
G_{N}}E\bigl(\bigl[\rho\bigl(\varepsilon_{\mathbf{i}}+e_{\mathbf{i}}-t^{T}Z_{\mathbf{i}}\bigr)-
\rho(\varepsilon_{\mathbf{i}}+e_{\mathbf{i}})\bigr]|X_{\mathbf
{i}},U_{\mathbf{i}}\bigr)K\biggl(\frac{U_{\mathbf{i}}
-u_{0}}{h_{\mathbf{n}}}\biggr).
\]
Let $\Upsilon_{\mathbf{n}}(t)=\sum_{\mathbf{i}\in
G_{N}}t^{T}Z_{\mathbf{i}}\psi(\varepsilon_{\mathbf{i}})
K(\frac{U_{\mathbf{i}}-u_{0}}{h_{\mathbf{n}}})$ and
$R_{\mathbf{n}}(t)=S_{\mathbf{n}}(t)-\Gamma_{\mathbf{n}}(t)+\Upsilon
_{\mathbf{n}}(t)$.
Then
%
\begin{equation}\label{A.3}
S_{\mathbf{n}}(t)=\Gamma_{\mathbf{n}}(t)-\Upsilon
_{\mathbf{n}}(t)+R_{\mathbf{n}}(t).
\end{equation}

We first present several lemmas that are necessary to prove the theorems.
%
\begin{Lem}\label{lem1}
Under the Assumptions \ref{asum1} and \ref{asum3}--\ref{asum8}, if
$n_{k}h_{\mathbf{n}}^{\delta
k/[a(2+\delta)]}>1$, then for any
fixed $t$, as $\mathbf{n}\rightarrow\infty$, it holds that
\[
R_{\mathbf{n}}(t)=\mathrm{o}_{p}(1).
\]
\end{Lem}
%
\begin{Lem}\label{lem2}
Assume that Assumptions \ref{asum1}, \ref{asum3} and \ref{asum5}--\ref{asum8} hold
and $\Phi(u)$ is continuous in some neighborhood of $u_{0}$. If
$n_{k}h_{\mathbf{n}}^{\delta k/[a(2+\delta)]}>1$, then, as
$\mathbf{n}\rightarrow\infty$, it holds that
\[
\Gamma_{\mathbf{n}}(t) =\tfrac{1}{2}f(u_{0})t^{T}\bigl(\Delta\otimes\Phi(u_{0})\bigr)t-
\tfrac{1}{2}\tilde{\mathbf{n}}^{1/2}h_{\mathbf{n}}^{(k+4)/2}f(u_{0})t^{T}
\Lambda\otimes\bigl(\Phi(u_{0})\zeta(u_{0})\bigr)+\mathrm{o}_{p}(1),
\]
where
$\Delta=\diag(\int_{\tilde{M}}K(u)\,\mathrm{d}u,\int_{\tilde{M}}uu^{T}K(u)\,\mathrm{d}u)$
and $\Lambda=(1,O^{T})^{T}$, where $O$ is a $k\times1$ vector
with entries zero.
\end{Lem}
%
\begin{Lem}\label{lem3}
Assume that Assumptions \ref{asum1}, \ref{asum3} and \ref{asum5}--\ref{asum8} hold
and $\Sigma(u)$ is continuous in some neighborhood
of $u_{0}$. If
$n_{k}h_{\mathbf{n}}^{\delta k/[a(2+\delta)]}>1$, then, as
$\mathbf{n}\rightarrow\infty$, it holds that
\[
D\bigl(\Upsilon_{\mathbf{n}}(t)\bigr) =f(u_{0})t^{T}\bigl(\tilde{\Delta}\otimes\Sigma
(u_{0})\bigr)t+\mathrm{o}(1),
\]
where
$\tilde{\Delta}=\diag(\int_{\tilde{M}}K^{2}(u)\,\mathrm{d}u,\int_{\tilde
{M}}uu^{T}K^{2}(u)\,\mathrm{d}u)$.
\end{Lem}
%
\begin{Lem}\label{lem4} Let
$K_{h}(U_{\mathbf{i}})=K(\frac{U_{\mathbf{i}}-u_{0}}{h_{\mathbf{n}}})$
and
\begin{eqnarray*}
A_{\mathbf{n}}(t)&=&\sum_{\mathbf{i}\in
G_{N}}\psi\bigl(\varepsilon_{\mathbf{i}}+e_{\mathbf{i}}-t^{T}Z_{\mathbf{i}}\bigr)
t^{T}\bigl[\bigl(1,h_{\mathbf{n}}^{-1}(U_{\mathbf{i}}-u_{0})^{T}\bigr)^{T}\otimes
(\hat{X}_{\mathbf{i}}-X_{\mathbf{i}})\bigr]K_{h}(U_{\mathbf{i}}),
\\
B_{\mathbf{n}}(t)&=&\sum_{\mathbf{i}\in
G_{N}}\psi\bigl(\varepsilon_{\mathbf{i}}+e_{\mathbf{i}}-t^{T}Z_{\mathbf{i}}\bigr)
t^{T}\bigl[\bigl(0,h_{\mathbf{n}}^{-1}(\hat{U}_{\mathbf{i}}-U_{\mathbf
{i}})^{T}\bigr)^{T}\otimes
X_{\mathbf{i}}\bigr]K_{h}(U_{\mathbf{i}}).
\end{eqnarray*}
Under the assumptions of Theorem~\ref{thm4.1}, for any fixed $t$, as
$\mathbf{n}\rightarrow\infty$, it holds that
\[
\bigl(\tilde{\mathbf{n}}h_{\mathbf{n}}^{k}\bigr)^{-1/2}A_{\mathbf
{n}}(t)=\mathrm{o}_{p}(1),\qquad  \bigl(\tilde{\mathbf{n}}h_{\mathbf
{n}}^{k}\bigr)^{-1/2}B_{\mathbf{n}}(t)=\mathrm{o}_{p}(1).
\]
\end{Lem}
%
\begin{Lem}\label{lem5}
Under the assumptions of
Theorem~\ref{thm4.1}, for any fixed $t$, as
$\mathbf{n}\rightarrow\infty$, it holds that
\[
\sum_{\mathbf{i}\in
G_{N}}\bigl[\rho\bigl(\varepsilon_{\mathbf{i}}+e_{\mathbf{i}}-t^{T}Z_{\mathbf{i}}\bigr)-
\rho(\varepsilon_{\mathbf{i}}+e_{\mathbf{i}})\bigr]\bigl[K_{h}(\hat{U}_{\mathbf{i}})
-K_{h}(U_{\mathbf{i}})\bigr]=\mathrm{o}_{p}(1).
\]
\end{Lem}
\begin{pf*}{Proof of Theorem~\ref{thm2.1}}
Let $\Gamma(t)
=\frac{1}{2}f(u_{0})t^{T}(\Delta\otimes\Phi(u_{0}))t$. By
Lemmas \ref{lem1} and \ref{lem2} and (\ref{A.3}), for fixed $t$
we have
%
\begin{equation}\label{A.4}
S_{\mathbf{n}}(t)=\Gamma(t)-\tfrac{1}{2}\tilde{\mathbf
{n}}^{1/2}h_{\mathbf{n}}^{(k+4)/2}f(u_{0})t^{T}
\Lambda\otimes
\bigl(\Phi(u_{0})\zeta(u_{0})\bigr)-\Upsilon_{\mathbf{n}}(t)+\tilde{R}_{\mathbf{n}}(t),
\end{equation}
where $\tilde{R}_{\mathbf{n}}(t)=R_{\mathbf{n}}(t)+\mathrm{o}_{p}(1)=\mathrm{o}_{p}(1)$,
and hence
\[
S_{\mathbf{n}}(t)+\tfrac{1}{2}\tilde{\mathbf{n}}^{1/2}h_{\mathbf
{n}}^{(k+4)/2}f(u_{0})t^{T}
\Lambda\otimes
\bigl(\Phi(u_{0})\zeta(u_{0})\bigr)+\Upsilon_{\mathbf{n}}(t)=\Gamma(t)+\tilde
{R}_{\mathbf{n}}(t).
\]
By Lemma~\ref{lem3}, $\Upsilon_{\mathbf{n}}(t)$ is bounded in
probability. Thus\DIFadd{,} the random convex
function $S_{\mathbf{n}}(t)+\tilde{\Gamma}_{\mathbf{n}}(t)+\Upsilon
_{\mathbf{n}}(t)$, for fixed $t$ converges
in probability to the function $\Gamma(t)$. According
to the convexity lemma \cite{s33}\DIFadd{,} we conclude that
for any
compact set $K$
%
\begin{equation}\label{A.5}
\sup_{t\in K}\bigl|\tilde{R}_{\mathbf{n}}(t)\bigr|=\mathrm{o}_{p}(1).
\end{equation}
Let $t^{\ast}=\frac{1}{2}\tilde{\mathbf{n}}^{1/2}h_{\mathbf
{n}}^{(k+4)/2}(\Delta^{-1}\Lambda)\otimes\zeta(u_{0})$
and
%
\begin{equation}\label{A.6}
\tilde{t}=t^{\ast}+\frac{1}{f(u_{0})}\bigl(\Delta^{-1}\otimes
\Phi^{-1}(u_{0})\bigr)
\sum_{\mathbf{i}\in G_{N}}Z_{\mathbf{i}}\psi(\varepsilon_{\mathbf
{i}})K\biggl(\frac{U_{\mathbf{i}}-u_{0}}{h_{\mathbf{n}}}\biggr).
\end{equation}
In the following\DIFadd{,} we will prove that for any sufficient small
$\epsilon>0$,
%
\begin{equation}\label{A.7}
P\bigl\{\|\hat{t}-\tilde{t}\|<\epsilon\bigr\}\rightarrow1.
\end{equation}
According to (\ref{A.1}) and Lemma~\ref{lem3} and using the convexity
of $\rho$,
to prove (\ref{A.7})\DIFadd{,} we need only to show that for any sufficient
large $L^{\ast}>0$,
%
\begin{equation}\label{A.8}
P\Bigl(\Bigl\{\inf_{\|t-\tilde{t}\|=\epsilon}\bigl(S_{\mathbf{n}}(t)-S_{\mathbf
{n}}(\tilde{t})\bigr)>0\Bigr\}\cap\bigl\{\|\tilde{t}\|\leq
L^{\ast}\bigr\}\Bigr)
\rightarrow1.
\end{equation}
By (\ref{A.4}) and (\ref{A.6}), we get
\[
\begin{array}{ll}
S_{\mathbf{n}}(t)=\frac{1}{2}f(u_{0})t^{T}\bigl(\Delta\otimes
\Phi(u_{0})\bigr)t-f(u_{0})t^{T}\bigl(\Delta\otimes
\Phi(u_{0})\bigr)\tilde{t}+\tilde{R}_{\mathbf{n}}(t).
\end{array}
\]
Since
\[
t^{T}\bigl(\Delta\otimes
\Phi(u_{0})\bigr)\tilde{t}=\tfrac{1}{2}\bigl[t^{T}\bigl(\Delta\otimes
\Phi(u_{0})\bigr)t+\tilde{t}^{T}\bigl(\Delta\otimes
\Phi(u_{0})\bigr)\tilde{t}-(t-\tilde{t})^{T}\bigl(\Delta\otimes
\Phi(u_{0})\bigr)(t-\tilde{t})\bigr].
\]
Hence\DIFadd{,}
\[
S_{\mathbf{n}}(t)=\tfrac{1}{2}f(u_{0})(t-\tilde{t})^{T}\bigl(\Delta\otimes
\Phi(u_{0})\bigr)(t-\tilde{t})-\tfrac{1}{2}f(u_{0})\tilde{t}^{T}\bigl(\Delta\otimes
\Phi(u_{0})\bigr)\tilde{t}+\tilde{R}_{\mathbf{n}}(t).
\]
Using the above, for $\tilde{t}$ satisfying that
$\|\tilde{t}\|\leq L^{\ast}$, it holds that
\[
S_{\mathbf{n}}(\tilde{t})=-\tfrac{1}{2}f(u_{0})\tilde{t}^{T}\bigl(\Delta
\otimes
\Phi(u_{0})\bigr)\tilde{t}+\tilde{R}_{\mathbf{n}}(\tilde{t}).
\]
Note that $\|t-\tilde{t}\|=\epsilon$, we conclude that
\[
S_{\mathbf{n}}(t)-S_{\mathbf{n}}(\tilde{t}) \geq
\tfrac{1}{2}f(u_{0})\lambda_{\min,\Delta}\lambda_{\min}(u_{0})\epsilon^{2}-
2\sup_{\|t\|\leq L^{\ast}+\epsilon} \bigl|\tilde{R}_{\mathbf{n}}(t)\bigr|,
\]
where $\lambda_{\min,\Delta}$ and $\lambda_{\min}(u_{0})$ are the
minimum eigenvalue of $\Delta$ and $\Phi(u_{0})$ respectively.
Therefore\DIFadd{,} (\ref{A.8}) follows from (\ref{A.5}) and the above.
Consequently\DIFadd{,} (\ref{A.7}) holds. Under assumptions of
Theorem~\ref{thm2.1}, using arguments similar to those used in the
proof of Lemma~3.1 of \cite{s14} and Lemma~\ref{lem3}, we can
show that
\[
\sum_{\mathbf{i}\in
G_{N}}Z_{\mathbf{i}}\psi(\varepsilon_{\mathbf{i}})K\biggl(\frac{U_{\mathbf
{i}}-u_{0}}{h_{\mathbf{n}}}\biggr)\rightarrow_{d}N\bigl(0,f(u_{0})\tilde{\Delta
}\otimes
\Sigma(u_{0})\bigr).
\]
Now the conclusion of Theorem~\ref{thm2.1} follows from
(\ref{A.2}), (\ref{A.7}), (\ref{A.6}) and the above, and the proof
of Theorem~\ref{thm2.1} is finished.
\end{pf*}
\begin{pf*}{Proof of Theorem~\ref{thm2.2}}
 Let
$\Pi_{0}(u_{0})=\Phi(u_{0})\bar{\zeta}(u_{0})$,
$\Pi_{l}(u_{0})=\Phi(u_{0})\bar{\zeta}^{(l)}(u_{0})$, $l=1,\ldots, k$,
$\Pi(u_{0})=(\Pi_{0}(u_{0})^{T},\Pi_{1}(u_{0})^{T},\ldots,
\Pi_{k}(u_{0})^{T})^{T}$ and
\begin{eqnarray*}
t^{\ast}&=&\frac{1}{2}\tilde{\mathbf{n}}^{1/2}h_{\mathbf
{n}}^{(k+4)/2}\bigl(\Delta_{c}^{-1}\otimes
\Phi^{-1}(u_{0})\bigr)\Pi(u_{0}),
\\
\tilde{t}&=&t^{\ast}+\frac{1}{f(u_{*})}\bigl(\Delta_{c}^{-1}\otimes\Phi^{-1}(u_{0})\bigr)
\sum_{\mathbf{i}\in G_{N}}Z_{\mathbf{i}}\psi(\varepsilon_{\mathbf
{i}})K\biggl(\frac{U_{\mathbf{i}}-u_{h}}{h_{\mathbf{n}}}\biggr).
\end{eqnarray*}
Using the arguments similar to those in the proof of Theorem~\ref{thm2.1}, we
can finish the proof of Theorem~\ref{thm2.2}.
\end{pf*}
\begin{pf*}{Proof of Theorem~\ref{thm4.1}}
Recall that $N = 2$ has been assumed
throughout this
section. Following \cite{s18}, $Y(s)$, $X(s)$ and $U(s)$ satisfy that
$\sup_{s\in
[0,1]^{2}}|\hat{\alpha}_{\mathbf{Y}}(s)-\alpha_{\mathbf
{Y}}(s)|=\mathrm{O}_{p}(\epsilon_{\mathbf{n}})$,
$\sup_{s\in
[0,1]^{2}}\|\hat{\alpha}_{\mathbf{X}}(s)-\alpha_{\mathbf{X}}(s)\|
=\mathrm{O}_{p}(\epsilon_{\mathbf{n}})$
and $ \sup_{s\in
[0,1]^{2}}\|\hat{\alpha}_{\mathbf{U}}(s)-\alpha_{\mathbf{U}}(s)\|
=\mathrm{O}_{p}(\epsilon_{\mathbf{n}})$
with
$\epsilon_{\mathbf{n}}=(\ln\tilde{\mathbf{n}}/(\tilde{\mathbf
{n}}g_{\mathbf{n}}^{2}))^{1/2}
+g_{\mathbf{n}}^{m}=:\epsilon_{\mathbf{n}}^{(1)}+\epsilon_{\mathbf{n}}^{(2)}$,
where $\epsilon_{\mathbf{n}}^{(1)}$ is obtained as in the proof of
Theorem~2 of \cite{s18} under Assumptions (B1)--(B3), (B5) and
\ref{asum8}, while $\epsilon_{\mathbf{n}}^{(2)}$ readily follows
from Assumptions (B3) and (B5). Therefore, we have
%
\begin{equation}\label{A.9}
\max_{\mathbf{i}}|\hat{Y}_{\mathbf{i}}-Y_{\mathbf{i}}|=\mathrm{O}_{p}(\epsilon
_{\mathbf{n}}),\qquad
\max_{\mathbf{i}}\|\hat{X}_{\mathbf{i}}-X_{\mathbf{i}}\|=\mathrm{O}_{p}(\epsilon
_{\mathbf{n}}),\qquad
\max_{\mathbf{i}}\|\hat{U}_{\mathbf{i}}-U_{\mathbf{i}}\|=\mathrm{O}_{p}(\epsilon
_{\mathbf{n}}).
\end{equation}
Let
$\hat{\varepsilon}_{\mathbf{i}}=\hat{Y}_{\mathbf{i}}-\hat{X}_{\mathbf
{i}}^{T}\beta(\hat{U}_{\mathbf{i}})$,
$\hat{e}_{\mathbf{i}}=\frac{1}{2}\beta^{\ast}(\bar{\xi}_{\mathbf{i}})^{T}
\hat{X}_{\mathbf{i}}$,
$\hat{Z}_{\mathbf{i}}=(\widetilde{\mathbf{n}}h_{\mathbf
{n}}^{k})^{-1/2}(1,h^{-1}(\hat{U}_{\mathbf{i}}-u_{0})^{T})^{T}\otimes
\hat{X}_{\mathbf{i}}$ and
$\check{t}=(\widetilde{\mathbf{n}}h_{\mathbf{n}}^{k})^{1/2}((\check
{a}-\beta(u_{0}))^{T},
h_{\mathbf{n}}(\check{b}-\dot{\beta}(u_{0}))^{T})^{T}$.
Then
\[
\check{t}=\mathop{\Argmin}_{t}\sum_{\mathbf{i}\in
G_{N}}\bigl[\rho\bigl(\hat{\varepsilon}_{\mathbf{i}}+\hat{e}_{\mathbf
{i}}-t^{T}\hat{Z}_{\mathbf{i}}\bigr)-
\rho(\hat{\varepsilon}_{\mathbf{i}}+\hat{e}_{\mathbf{i}})\bigr]K_{h}(\hat
{U}_{\mathbf{i}}),
\]
where
$K_{h}(\hat{U}_{\mathbf{i}})=K((\hat{U}_{\mathbf{i}}-u_{0})/h_{\mathbf{n}})$.
Let $ \hat{S}_{\mathbf{n}}(t)=\sum_{\mathbf{i}\in
G_{N}}[\rho(\hat{\varepsilon}_{\mathbf{i}}+\hat{e}_{\mathbf
{i}}-t^{T}\hat{Z}_{\mathbf{i}})-
\rho(\hat{\varepsilon}_{\mathbf{i}}+\hat{e}_{\mathbf{i}})]K_{h}(\hat
{U}_{\mathbf{i}})$.
According to (\ref{A.4}) and the proof of Theorem~\ref{thm2.1}, to
finish the proof of Theorem~\ref{thm4.1}, we need only show that
for fixed $t$, it holds that
%
\begin{equation}\label{A.10}
\hat{S}_{\mathbf{n}}(t)-S_{\mathbf{n}}(t)=\mathrm{o}_{p}(1).
\end{equation}
Let
$\theta_{\mathbf{i}}=\varepsilon_{\mathbf{i}}+e_{\mathbf{i}},
\hat{\theta}_{\mathbf{i}}=\hat{\varepsilon}_{\mathbf{i}}+\hat
{e}_{\mathbf{i}}$,
$V_{\mathbf{n}1}=\sum_{\mathbf{i}\in
G_{N}}[(\rho(\hat{\theta}_{\mathbf{i}}-t^{T}\hat{Z}_{\mathbf{i}})-
\rho(\hat{\theta}_{\mathbf{i}}))-(\rho(\theta_{\mathbf
{i}}-t^{T}Z_{\mathbf{i}})-
\rho(\theta_{\mathbf{i}}))]K_{h}(\hat{U}_{\mathbf{i}})$,
and $V_{\mathbf{n}2}=\sum_{\mathbf{i}\in
G_{N}}[\rho(\theta_{\mathbf{i}}-t^{T}Z_{\mathbf{i}})-
\rho(\theta_{\mathbf{i}})][K_{h}(\hat{U}_{\mathbf{i}})-K_{h}(U_{\mathbf{i}})]$.
Then
%
\begin{equation}\label{A.11}
\hat{S}_{\mathbf{n}}(t)-S_{\mathbf{n}}(t)=V_{\mathbf{n}1}+V_{\mathbf{n}2}.
\end{equation}
Let $V_{\mathbf{i}1}=|\psi(\hat{\theta}_{\mathbf{i}}-t^{T}\hat
{Z}_{\mathbf{i}})
-\psi(\theta_{\mathbf{i}}-t^{T}Z_{\mathbf{i}})|
\cdot
|(\hat{\theta}_{\mathbf{i}}-\theta_{\mathbf{i}})-t^{T}(\hat{Z}_{\mathbf
{i}}-Z_{\mathbf{i}})|$,
$V_{\mathbf{i}2}=|\psi(\hat{\theta}_{\mathbf{i}})-\psi(\theta_{\mathbf{i}})|
\cdot|(\hat{\theta}_{\mathbf{i}}-\theta_{\mathbf{i}})|$ and
$V_{\mathbf{i}3}=|\psi(\theta_{\mathbf{i}}-t^{T}Z_{\mathbf{i}})-\psi
(\theta_{\mathbf{i}})|
\cdot|(\hat{\theta}_{\mathbf{i}}-\theta_{\mathbf{i}})|$. By the
convexity of $\rho(\cdot)$, it holds that
\[
\bigl|\rho\bigl(\hat{\theta}_{\mathbf{i}}-t^{T}\hat{Z}_{\mathbf{i}}\bigr)-
\rho\bigl(\theta_{\mathbf{i}}-t^{T}Z_{\mathbf{i}}\bigr)-\psi\bigl(\theta_{\mathbf
{i}}-t^{T}Z_{\mathbf{i}}\bigr)
\bigl[(\hat{\theta}_{\mathbf{i}}-\theta_{\mathbf{i}})-t^{T}(\hat{Z}_{\mathbf
{i}}-Z_{\mathbf{i}})\bigr]\bigr|
\leq V_{\mathbf{i}1}
\]
and $|\rho(\hat{\theta}_{\mathbf{i}})-
\rho(\theta_{\mathbf{i}})-\psi(\theta_{\mathbf{i}})
(\hat{\theta}_{\mathbf{i}}-\theta_{\mathbf{i}})|
\leq V_{\mathbf{i}2}$.
Hence
%
\begin{equation}\label{A.12}
V_{\mathbf{n}1}\leq\sum_{\mathbf{i}\in
G_{N}}(V_{\mathbf{i}1}+V_{\mathbf{i}2}+V_{\mathbf{i}3})K_{h}(\hat
{U}_{\mathbf{i}})+|V_{\mathbf{n}3}|,
\end{equation}
where
%
\begin{equation}\label{A.13}
V_{\mathbf{n}3}=\sum_{\mathbf{i}\in
G_{N}}\psi\bigl(\theta_{\mathbf{i}}-t^{T}Z_{\mathbf{i}}\bigr)t^{T}(\hat
{Z}_{\mathbf{i}}
-Z_{\mathbf{i}})K_{h}(\hat{U}_{\mathbf{i}}).
\end{equation}
Since
$\theta_{\mathbf{i}}=\varepsilon_{\mathbf{i}}+e_{\mathbf{i}}=Y_{\mathbf
{i}}-X_{\mathbf{i}}^{T}\beta(u_{0})
-\sum_{r=1}^{d}(U_{\mathbf{i}}-u_{0})^{T}\dot{\beta
}_{r}(u_{0})X_{\mathbf{i}r}$
and
$\hat{\theta}_{\mathbf{i}}=\hat{Y}_{\mathbf{i}}-\hat{X}_{\mathbf
{i}}^{T}\beta(u_{0})
-\sum_{r=1}^{d}(\hat{U}_{\mathbf{i}}-u_{0})^{T}\dot{\beta
}_{r}(u_{0})\hat{X}_{\mathbf{i}r}$,
by (\ref{A.9}) and Assumption~\ref{asum6'}, it is easy to prove that
$\max_{\mathbf{i}}|\hat{\theta}_{\mathbf{i}}-\theta_{\mathbf
{i}}|=\mathrm{O}_{p}(\epsilon_{\mathbf{n}})$.
On the other hand,
\[
t^{T}(\hat{Z}_{\mathbf{i}}-Z_{\mathbf{i}})=\bigl(\tilde{\mathbf{n}}h_{\mathbf
{n}}^{k}\bigr)^{-1/2}
\Biggl[\Biggl(t_{0}+\sum_{l=1}^{k}\frac{\hat{U}_{\mathbf{i}l}-u_{0l}}{h_{\mathbf
{n}}}t_{l}\Biggr)^{T}(\hat{X}_{\mathbf{i}}-
X_{\mathbf{i}})+\Biggl(\sum_{l=1}^{k}\frac{\hat{U}_{\mathbf{i}l}-U_{\mathbf
{i}l}}{h_{\mathbf{n}}}t_{l}\Biggr)^{T}X_{\mathbf{i}}\Biggr].
\]
By (\ref{A.9}) and Assumption~\ref{asum6'}, we have
$\max_{\mathbf{i}}|t^{T}(\hat{Z}_{\mathbf{i}}-Z_{\mathbf{i}})|
=\mathrm{O}_{p}((\tilde{\mathbf{n}}h_{\mathbf{n}}^{k})^{-1/2}h_{\mathbf
{n}}^{-1}\epsilon_{\mathbf{n}})$.
Hence
%
\begin{equation}\label{A.14}
\max_{\mathbf{i}}\bigl(|\hat{\theta}_{\mathbf{i}}-\theta
_{\mathbf{i}}|+\bigl|t^{T}(\hat{Z}_{\mathbf{i}}-Z_{\mathbf{i}})\bigr|\bigr)
=\mathrm{O}_{p}\bigl(\epsilon_{\mathbf{n}}+\bigl(\tilde{\mathbf{n}}h_{\mathbf
{n}}^{k}\bigr)^{-1/2}h_{\mathbf{n}}^{-1}\epsilon_{\mathbf{n}}\bigr)
=\mathrm{O}_{p}(\tilde{\epsilon}_{\mathbf{n}}),
\end{equation}
where $\tilde{\epsilon}_{\mathbf{n}}=\epsilon_{\mathbf{n}}+
(\tilde{\mathbf{n}}h_{\mathbf{n}}^{k})^{-1/2}h_{\mathbf{n}}^{-1}\epsilon
_{\mathbf{n}}$.
By Assumption~\ref{asum7'}, we get
%
\begin{equation}\label{A.15}
K_{h}(\hat{U}_{\mathbf{i}})=K_{h}(U_{\mathbf{i}})
+h_{\mathbf{n}}^{-1}(\hat{U}_{\mathbf{i}}-U_{\mathbf{i}})^{T}\dot
{K}_{h}(U_{\mathbf{i}})\bigl[1+\mathrm{o}_{p}(1)\bigr]
=K_{h}(U_{\mathbf{i}})+\mathrm{o}_{p}(1).
\end{equation}
Therefore\DIFadd{, }
%
\begin{eqnarray}\label{A.16}
\sum_{\mathbf{i}\in
G_{N}}V_{\mathbf{i}1}K_{h}(\hat{U}_{\mathbf{i}}) & =&\bigl[1+\mathrm{o}_{p}(1)\bigr]\sum
_{\mathbf{i}\in
G_{N}}V_{\mathbf{i}1}K_{h}(U_{\mathbf{i}}) \nonumber\\[-8pt]\\[-8pt]
&
=&\mathrm{O}_{p}(\tilde{\epsilon}_{\mathbf{n}})
\sum_{\mathbf{i}\in
G_{N}}\bigl|\psi\bigl(\hat{\theta}_{\mathbf{i}}-t^{T}\hat{Z}_{\mathbf{i}}\bigr)
-\psi\bigl(\theta_{\mathbf{i}}-t^{T}Z_{\mathbf{i}}\bigr)\bigr|K_{h}(U_{\mathbf{i}}).\nonumber
\end{eqnarray}
According to (\ref{A.14}), we can assume that, with probability
arbitrarily close to one,
$\max_{\mathbf{i}}(|\hat{\theta}_{\mathbf{i}}-\theta_{\mathbf
{i}}|+|t^{T}(\hat{Z}_{\mathbf{i}}-Z_{\mathbf{i}})|)
\leq C\tilde{\epsilon}_{\mathbf{n}}$ for some $C$ and $\mathbf{n}$
sufficiently large. Then by Assumption~\ref{asum4'}, it holds that
$\sum_{\mathbf{i}\in
G_{N}}E(E(|\psi(\hat{\theta}_{\mathbf{i}}-t^{T}\hat{Z}_{\mathbf{i}})
-\psi(\theta_{\mathbf{i}}-t^{T}Z_{\mathbf{i}})||U_{\mathbf
{i}})K_{h}(U_{\mathbf{i}}))
=\mathrm{O}(\tilde{\mathbf{n}}h_{\mathbf{n}}^{k}\tilde{\epsilon}_{\mathbf{n}})$.
Therefore, by\vspace*{1pt} Assumption~\ref{asum5'}, it holds that
%
\begin{equation}\label{A.17}
\sum_{\mathbf{i}\in
G_{N}}V_{\mathbf{i}1}K_{h}(\hat{U}_{\mathbf{i}})
=\mathrm{O}_{p}\bigl(\tilde{\mathbf{n}}h_{\mathbf{n}}^{k}\tilde{\epsilon}_{\mathbf{n}}^{2}\bigr)
=\mathrm{O}_{p}\bigl(\tilde{\mathbf{n}}h_{\mathbf{n}}^{k}\epsilon_{\mathbf
{n}}^{2}+h_{\mathbf{n}}^{-2}\epsilon_{\mathbf{n}}^{2}\bigr)
=\mathrm{o}_{p}(1).
\end{equation}
Similarly
%
\begin{equation}\label{A.18}
\sum_{\mathbf{i}\in
G_{N}}V_{\mathbf{i}2}K_{h}(\hat{U}_{\mathbf{i}})
=\mathrm{O}_{p}\bigl(\tilde{\mathbf{n}}h_{\mathbf{n}}^{k}\epsilon_{\mathbf{n}}^{2}\bigr)
=\mathrm{o}_{p}(1)
\end{equation}
and
%
\begin{equation}\label{A.19}
 \sum_{\mathbf{i}\in
G_{N}}V_{\mathbf{i}3}K_{h}(\hat{U}_{\mathbf{i}})
=\mathrm{O}_{p}\biggl(\tilde{\mathbf{n}}h_{\mathbf{n}}^{k}\max_{\mathbf
{i}}\bigl|t^{T}Z_{\mathbf{i}}\bigr|\epsilon_{\mathbf{n}}\biggr)
=\mathrm{O}_{p}\bigl(\bigl(\tilde{\mathbf{n}}h_{\mathbf{n}}^{k}\bigr)^{1/2}\epsilon_{\mathbf
{n}}\bigr)=\mathrm{o}_{p}(1).
\end{equation}
By (\ref{A.13}), (\ref{A.15}) and Lemma~\ref{lem4}, we obtain
%
\begin{eqnarray}\label{A.20}
V_{\mathbf{n}3} & =&\bigl[1+\mathrm{o}_{p}(1)\bigr]\sum_{\mathbf{i}\in
G_{N}}\psi\bigl(\theta_{\mathbf{i}}-t^{T}Z_{\mathbf{i}}\bigr)t^{T}(\hat
{Z}_{\mathbf{i}}
-Z_{\mathbf{i}})K_{h}(U_{\mathbf{i}})\nonumber \\[-8pt]\\[-8pt]
& =&\bigl[1+\mathrm{o}_{p}(1)\bigr]\bigl[\bigl(\tilde{\mathbf{n}}h_{\mathbf{n}}^{k}\bigr)^{-1/2}A_{\mathbf{n}}(t)+
\bigl(\tilde{\mathbf{n}}h_{\mathbf{n}}^{k}\bigr)^{-1/2}B_{\mathbf{n}}(t)\bigr]=\mathrm{o}_{p}(1).\nonumber
\end{eqnarray}
Combining (\ref{A.12}) and
(\ref{A.17})--(\ref{A.20}), we conclude that
$V_{\mathbf{n}1}=\mathrm{o}_{p}(1)$. By Lemma~\ref{lem5}, it holds that
$V_{\mathbf{n}2}=\mathrm{o}_{p}(1)$. Therefore, by (\ref{A.11}),
(\ref{A.10}) holds and the proof of Theorem~\ref{thm4.1} is
finished.
\end{pf*}
\begin{pf*}{Proof of Theorem~\ref{thm4.2}} The proof of Theorem~\ref{thm4.2}
is similar
to that of Theorem~\ref{thm4.1} except proof of (\ref{A.17}). Let
$\vartheta_{\mathbf{i}}=\hat{\theta}_{\mathbf{i}}-
\theta_{\mathbf{i}}-t^{T}(\hat{Z}_{\mathbf{i}}-Z_{\mathbf{i}})$. Since
$\psi(z)=2\tau
I(z>0)+2(\tau-1)I(z<0)$, it holds that
\[
\bigl|\psi\bigl(\hat{\theta}_{\mathbf{i}}-t^{T}\hat{Z}_{\mathbf{i}}\bigr)
-\psi\bigl(\theta_{\mathbf{i}}-t^{T}Z_{\mathbf{i}}\bigr)\bigr| \leq
2I_{\{|\theta_{\mathbf{i}}-t^{T}Z_{\mathbf{i}}|\leq
|\vartheta_{\mathbf{i}}|\}} \leq
2I_{\{|\varepsilon_{\mathbf{i}}|\leq
|e_{\mathbf{i}}|+|t^{T}Z_{\mathbf{i}}|+|\vartheta_{\mathbf{i}}|\}}.
\]
By Assumptions \ref{asum5'} and \ref{asum6'} and (\ref{A.14}), we have
$\max_{\mathbf{i}}(|e_{\mathbf{i}}|+|t^{T}Z_{\mathbf{i}}|+|\vartheta
_{\mathbf{i}}|)
=\mathrm{O}_{p}((\tilde{\mathbf{n}}h_{\mathbf{n}}^{k})^{-1/2}+\tilde{\epsilon
}_{\mathbf{n}})$.
Thus we can assume that, with probability arbitrarily close to
one,\break
$\max_{\mathbf{i}}(|e_{\mathbf{i}}|+|t^{T}Z_{\mathbf{i}}|+|\vartheta
_{\mathbf{i}}|)
\leq
C((\tilde{\mathbf{n}}h_{\mathbf{n}}^{k})^{-1/2}+\tilde{\epsilon
}_{\mathbf{n}})$
for some $C$ and $\mathbf{n}$ sufficiently large. By Assumption~\ref{asumQ}
and the fact that
$f_{\varepsilon}(0|X_{\mathbf{i}},U_{\mathbf{i}})\leq C$ for some
$C>0$, we get that $EI_{\{|\varepsilon_{\mathbf{i}}|\leq
C((\tilde{\mathbf{n}}h_{\mathbf{n}}^{k})^{-1/2}+\tilde{\epsilon
}_{\mathbf{n}})\}}K_{h}(U_{\mathbf{i}})
=\mathrm{O}(((\tilde{\mathbf{n}}h_{\mathbf{n}}^{k})^{-1/2}+\tilde{\epsilon
}_{\mathbf{n}})h_{\mathbf{n}}^{k})$.
Therefore
%
\begin{eqnarray}\label{A.21}
&&\sum_{\mathbf{i}\in
G_{N}}\bigl|\psi\bigl(\hat{\theta}_{\mathbf{i}}-t^{T}\hat{Z}_{\mathbf{i}}\bigr)
-\psi\bigl(\theta_{\mathbf{i}}-t^{T}Z_{\mathbf{i}}\bigr)\bigr|K_{h}(U_{\mathbf{i}})\nonumber\\[-8pt]\\[-8pt]
&&\quad \leq\bigl[1+\mathrm{o}_{p}(1)\bigr]\sum_{\mathbf{i}\in
G_{N}}EI_{\{|\varepsilon_{\mathbf{i}}|\leq
C((\tilde{\mathbf{n}}h_{\mathbf{n}}^{k})^{-1/2}+\tilde{\epsilon
}_{\mathbf{n}})\}}K_{h}(U_{\mathbf{i}})
=\mathrm{O}\bigl(\bigl(\tilde{\mathbf{n}}h_{\mathbf{n}}^{k}\bigr)^{1/2}+\tilde{\mathbf
{n}}h_{\mathbf{n}}^{k}\tilde{\epsilon}_{\mathbf{n}}\bigr).\nonumber
\end{eqnarray}
Hence by (\ref{A.16}), (\ref{A.21}) and Assumption~\ref{asum5'}, we obtain
\[
\sum_{\mathbf{i}\in
G_{N}}V_{\mathbf{i}1}K_{h}(\hat{U}_{\mathbf{i}})
=\mathrm{O}\bigl(\bigl(\tilde{\mathbf{n}}h_{\mathbf{n}}^{k}\bigr)^{1/2}\tilde{\epsilon
}_{\mathbf{n}}
+\tilde{\mathbf{n}}h_{\mathbf{n}}^{k}\tilde{\epsilon}_{\mathbf{n}}^{2}\bigr)
=\mathrm{o}_{p}(1).
\]
Therefore\DIFadd{, } (\ref{A.17}) holds and the proof of Theorem~\ref{thm4.2}
is finished.
\end{pf*}
\begin{pf*}{Proof of Theorem~\ref{thm4.3}} The proof of
Theorem~\ref{thm4.3} can be done similarly as in that for
Theorem~\ref{thm4.2}, and the detail is omitted.
\end{pf*}
\end{appendix}

\section*{Acknowledgements}
We would first of all express our gratitude to both anonymous
referees as well as the chief editor, Prof. Richard A. Davis, and
an associate editor for their valuable comments and suggestions,
which had greatly improved the early version of this paper. We
would also thank Zhenyu Jiang for the computational helps in
preparing the real data example in Section~\ref{sec5}. This research was
supported by the Australian Research Council's Discovery Project
Grant DP0984686 and Future Fellowships Grant FT100100109, which
are acknowledged. Tang's research was also partially supported by
National Natural Science Foundation of China (Grant 11071120).
\begin{supplement}
\sname{Supplement}
\stitle{Proofs of the lemmas in \hyperref[Aproofs]{Appendix}}
\slink[doi]{10.3150/12-BEJ480SUPP} 
\sdatatype{.pdf}
\sfilename{BEJ480\_supp.pdf}
\sdescription{We collect the proofs for the necessary lemmas used
in the above in this supplementary material \cite{s31a}.}
\end{supplement}


\printhistory

\end{document}